\documentclass[12pt,onecolumn]{IEEEtran}
\usepackage[utf8]{inputenc}
\usepackage{amsmath}
\usepackage{amsfonts}
\usepackage{amssymb}
\usepackage{graphicx}
\usepackage{array}
\usepackage{multirow}
\usepackage[table]{xcolor}
\usepackage{listingsutf8}
\usepackage{geometry}
\usepackage{optidef}
\usepackage{hyperref}
\usepackage{float}
\usepackage{microtype}
\usepackage[authoryear]{natbib}
\usepackage{enumitem}
\usepackage{tikz}
\usepackage{algorithm}
\usepackage{algpseudocode}

\newtheorem{theorem}{Theorem}
\makeatletter

\@addtoreset{equation}{section}
\makeatother

\title{Comparative Evaluation of SDP, SOCP, and QC Convex Relaxations for Large-Scale Market-Based AC Optimal Power Flow}
\author{Ata Keskin \\ Technical University of Munich}

\begin{document}

\maketitle

\begin{abstract}
The alternating current optimal power flow (ACOPF) problem is central to modern power system operations, determining how electricity is generated and transmitted to maximize social welfare while respecting physical and operational constraints. However, the nonlinear and non-convex nature of AC power flow equations makes finding globally optimal solutions computationally intractable for large networks. Convex relaxations—including semidefinite programming (SDP), second-order cone programming (SOCP), and quadratic convex (QC) formulations—provide tractable alternatives that can yield provably optimal or near-optimal solutions under appropriate conditions. This paper presents a comprehensive comparative study of multiple ACOPF relaxations applied to market-based welfare maximization. We implement DCOPF, Shor's SDP relaxation (complex and real-valued forms), chordal SDP, Jabr's SOCP relaxation, and QC relaxations in a unified, solver-native framework using the MOSEK Fusion API, eliminating modeling overhead present in high-level frameworks such as CVXPY. To address the practical challenge of missing or overly conservative angle difference bounds required by QC relaxations, we employ quasi-Monte Carlo sampling with Sobol sequences to empirically estimate tighter bounds. We evaluate these relaxations on subnetworks of varying sizes derived from the ARPA-E dataset, systematically comparing solution quality, runtime, and memory consumption. Our results demonstrate the trade-offs between relaxation tightness and computational efficiency, providing practical guidance for selecting appropriate formulations based on network scale and solution requirements.
\end{abstract}

\section{Introduction}
The alternating current optimal power flow (ACOPF) problem, particularly in a market-based setting, is a central optimization challenge in modern power system operations, where the goal is to allocate generation and consumption to maximize social welfare while respecting physical and engineering constraints. It governs how electricity is generated and transmitted across networks, ensuring that buyers’ demands are met and sellers’ costs are accounted for under network constraints. Despite its importance, ACOPF is notoriously difficult to solve: the underlying power flow equations are nonlinear and non-convex, which makes finding globally optimal solutions computationally intractable for large networks.

Over the past two decades, convex relaxations have emerged as powerful tools for tackling this non-convexity. Semidefinite programming (SDP) relaxations, second-order cone programming (SOCP) formulations, and quadratic convex (QC) relaxations provide convex surrogates of ACOPF that yield provable upper bounds on the optimal welfare \citep{Lavaei2012ZeroDG, jabr2006, coffrin2016}. Under specific conditions, these relaxations turn out to be exact, meaning that a feasible ACOPF solution attaining the global optimum can be recovered from the relaxation \citep{Lavaei2012ZeroDG}. Among them, SDP-based methods offer the strongest theoretical guarantees, SOCP formulations scale more effectively to large networks \citep{jabr2006}, and QC relaxations strengthen SOCP models by incorporating convex envelopes for the nonlinear trigonometric structure of ACOPF \citep{coffrin2016}.

This paper makes three main contributions to the study of convex relaxations for the market-based formulation of ACOPF. First, we develop a unified framework that specifies and implements the DC approximation as a baseline alongside several widely studied convex relaxations: Shor's SDP relaxation (both complex and real-valued forms), the chordal SDP variant, Jabr's SOCP relaxation, and the QC relaxation. We began by prototyping each model in CVXPY for ease of development, then reimplemented them in the MOSEK Fusion API to reduce modeling overhead and improve computational efficiency. Along the way, we document the key modeling decisions behind each formulation, with particular attention to the chordal SDP relaxation \citep{fukuda2001, andersen2015, Madani2015}—covering how chordal extensions are constructed and how the full voltage matrix can be recovered from clique-based constraints.

Second, we address the practical challenge that QC relaxations require bounds on voltage angle differences, which are sometimes unavailable or specified in overly conservative form. To overcome this, we use quasi-Monte Carlo sampling using Sobol sequences \citep{sobol1967, Caflisch1998} to empirically estimate lower and upper bounds on angle differences for each transmission line. These empirically derived ranges, together with the corresponding bounds on trigonometric values, allow us to define non-trivial convex envelopes even when no reliable analytical bounds are provided. In this way, we approximate a tighter QC relaxation in settings where standard formulations would otherwise rely on suboptimal or missing angle-difference limits.

Finally, we present a comparative study of these relaxations on subnetworks derived from the ARPA-E dataset \citep{elbert2024arpa}. By sampling connected subgraphs of varying sizes that contain at least one generator and one load, we assess the scaling behavior of different relaxations with increasing network size. Our experiments compare not only the tightness of relaxations but also their runtime and memory demands, thereby highlighting practical tradeoffs relevant to large-scale power system optimization.

\section{Related Work}

The study of convex relaxations for the ACOPF problem has expanded significantly over the past two decades. A landmark contribution came from Lavaei and Low \citet{Lavaei2012ZeroDG}, who showed that a semidefinite programming (SDP) relaxation could often recover globally optimal solutions for standard IEEE test networks. Their results also identified structural conditions under which the relaxation is guaranteed to be exact. This discovery motivated substantial follow-up work into SDP relaxations of ACOPF (e.g., \citet{Madani2015, molzahn2017}), which clarified both their strengths and their limitations: while they can be very tight, the main drawback is their lack of scalability, since enforcing a positive semidefinite constraint on the full system matrix quickly becomes computationally demanding. To improve scalability, several works have exploited chordal decomposition, enforcing positive semidefiniteness only on smaller cliques rather than the full matrix (\citet{fukuda2001, waki2006, jabr2012, Madani2015, andersen2015}).

To address the scalability limitations of SDP methods, second-order cone programming (SOCP) relaxations, introduced by \citet{jabr2006}, provide a more compact alternative by focusing on line-flow variables. These “edge-based” relaxations are more compact and can scale to much larger networks than full SDP models, but they generally yield weaker bounds.

Building on this, \citet{coffrin2016} proposed quadratic convex (QC) relaxations, which aim to strengthen SOCP formulations by more carefully approximating the trigonometric structure of ACOPF. QC relaxations replace nonlinearities with convex envelopes, such as McCormick relaxations for bilinear terms and convex approximations for sine and cosine functions. Later work improved QC relaxations further with better envelope constructions and systematic bound-tightening procedures \citep{hijazi2017, coffrin2021}.

Moment relaxations generalize the basic SDP approach by applying polynomial optimization theory (e.g., \citet{waki2006, molzahn2019}). These methods can recover 
the exact global solution at sufficiently high relaxation orders \citep{nie2014optimality, josz2018lasserre}, but they are only tractable for very small systems. Hybrid methods have also been developed, such as using SDP or SOCP-based tightening in preprocessing, or adding cutting planes to refine relaxations iteratively (e.g., \citet{ghaddar2016, Kocuk2016}).

Several studies have implemented ACOPF relaxations using high-level frameworks such as MATPOWER \citep{matpower2011}, PYPOWER \citep{Lincoln2011}, and PowerModels.jl \citep{coffrin2018powermodelsjl}. These tools facilitate rapid prototyping of SDP, SOCP, and QC relaxations, but canonicalization and abstraction layers can create overhead that limits large-scale experiments.

While several studies have compared relaxation approaches, comprehensive evaluations that simultaneously assess solution quality, runtime, and memory consumption across networks of varying sizes remain limited. Early benchmarks focused primarily on small IEEE test cases or individual networks \citep{Lavaei2012ZeroDG, coffrin2016, hijazi2017}, while more recent work has extended to larger networks and subnetworks \citep{molzahn2017, coffrin2021}. However, most studies emphasize solution quality or theoretical tightness, with limited attention to computational efficiency metrics or the overhead introduced by high-level modeling frameworks. Furthermore, these evaluations are typically conducted on cost-minimization formulations rather than market-based welfare maximization objectives.

Market-based ACOPF, which maximizes social welfare instead of minimizing generation cost, is less explored. Foundational work includes \citet{bitar2012}, \citet{kirschen2004, kirschen2004_2} and \citet{moreira2022}. While prior work incorporates detailed constraints, network security, and strategic bidding, our formulation focuses on maximizing social welfare with basic bus-level constraints.

Our work addresses these gaps by providing a unified solver-native implementation in MOSEK Fusion that eliminates modeling overhead present in high-level frameworks, applying these methods to the market-based ACOPF formulation, and conducting a systematic comparative evaluation on ARPA-E subnetworks of varying sizes that analyzes solution quality, runtime, and memory consumption.

\section{Problem Formulation}

We study a market-based AC optimal power flow (ACOPF) problem in which electricity is traded between sellers (generators) and buyers (loads) located at network buses. 

Let $\mathcal{G} = (\mathcal{V}, \mathcal{E})$ denote the power network, where $\mathcal{V}$ is the set of buses and $\mathcal{E}$ is the set of transmission lines. For a node $v \in \mathcal{V}$, let $\mathcal{N}(v)$ denote the neighbors of $v$. The network is assumed to be meshed, meaning that transmission lines form cycles and provide multiple paths between buses, as opposed to radial (tree-structured) networks. For simplicity, the model does not consider transformers or shunt elements. Let $\mathcal{B}$ be the set of buyers and $\mathcal{S}$ the set of sellers. For a node $v$, let $\mathcal{B}_v$ be the set of buyers at that node and $\mathcal{S}_v$ be the set of sellers at that node. Time is discretized into periods indexed by $t \in \mathcal{T}$. Market participants submit bids and offers in blocks ($L_\mathcal{B}$ and $L_\mathcal{S}$ respectively), each specifying a quantity of power and an associated valuation (for buyers) or cost (for sellers). The total demand or supply at a node is the sum of the allocated block volumes. This block structure allows generation and consumption to be represented by monotone, piecewise-linear functions. 

We denote the valuation of buyer $b \in \mathcal{B}$ for block $l \in L_\mathcal{B}$ at time $t \in \mathcal{T}$ with $v_{btl}$ and the block size with $\sigma_{btl}$. Similarly, $c_{stl}$ is the cost that seller $s \in \mathcal{S}$ demands for block $l \in L_\mathcal{S}$ at time $t \in \mathcal{T}$ and the block size is denoted as $\sigma_{stl}$. In many real world scenarios, it costs money to just run a generator without any net power output. We model this no-load cost for a seller $s$ at time $t$ via the constant $C^\text{no-load}_{st}$. 

We define $p_{bt}$ and $q_{bt}$ as the active and reactive consumption of buyer $b$ at time $t$, and similarly $p_{st}$ and $q_{st}$ as the production of seller $s$ at time $t$. The variable $p_{btl}$ represents the consumption corresponding to block $l$ for buyer $b$, and $p_{stl}$ represents the production for block $l$ by seller $s$. Active and reactive power are bounded by $[\underline{p}_{bt}, 
\overline{p}_{bt}]$ and $[\underline{q}_{bt}, \overline{q}_{bt}]$ for buyers, and $[\underline{p}_{st}, \overline{p}_{st}]$ and $[\underline{q}_{st}, \overline{q}_{st}]$ for sellers.

Finally, for each transmission line $(v,w) \in \mathcal{E}$, we represent the real and reactive flows from bus $v$ to bus $w$ with the variables $p_{vwt}$ and $q_{vwt}$. We deliberately overload the letters $p$ and $q$, but the meaning is always clear from the indices, allowing us to keep the notation compact and readable.

The market-clearing problem maximizes the sum of buyer valuations minus seller costs, subject to the following constraints:
\begin{maxi!}
  {\substack{p_{btl},\,p_{stl},\\p_{bt},\,p_{st},\,u_{st}}}
  { \sum_{b \in \mathcal{B}} \sum_{t \in \mathcal{T}} \sum_{l \in L_\mathcal{B}} v_{btl}\, p_{btl}
    - \sum_{s \in \mathcal{S}} \sum_{t \in \mathcal{T}} \sum_{l \in L_\mathcal{S}} c_{stl}\, p_{stl}
    - \sum_{s \in \mathcal{S}} \sum_{t \in \mathcal{T}} C^{\text{no-load}}_{st}\, u_{st}}
  {}{}
  \addConstraint{\sum_{l \in L_\mathcal{B}} p_{btl}}{= p_{bt}}{\hspace{6em}\forall b \in \mathcal{B},\, t \in \mathcal{T}}
  \addConstraint{p_{btl}}{\in [0, \sigma_{btl}]}{\hspace{6em}\forall b \in \mathcal{B},\, t \in \mathcal{T}, \, l \in L_\mathcal{B}}
  \addConstraint{p_{bt}}{\in [\underline{p}_{bt}, \overline{p}_{bt}]}{\hspace{6em}\forall b \in \mathcal{B},\, t \in \mathcal{T}}
  \addConstraint{q_{bt}}{\in [\underline{q}_{bt}, \overline{q}_{bt}]}{\hspace{6em}\forall b \in \mathcal{B},\, t \in \mathcal{T}}
  \addConstraint{\sum_{l \in L_\mathcal{S}} p_{stl}}{= p_{st}}{\hspace{6em}\forall s \in \mathcal{S},\, t \in \mathcal{T}}
  \addConstraint{p_{stl}}{\in [0, \sigma_{stl} \cdot u_{st}]}{\hspace{6em}\forall s \in \mathcal{S},\, t \in \mathcal{T}, \, l \in L_\mathcal{S}}
  \addConstraint{p_{st}}{\in [\underline{p}_{st} \cdot u_{st}, \overline{p}_{st} \cdot u_{st}]}{\hspace{6em}\forall s \in \mathcal{S},\, t \in \mathcal{T}}
  \addConstraint{q_{st}}{\in [\underline{q}_{st} \cdot u_{st}, \overline{q}_{st} \cdot u_{st}]}{\hspace{6em}\forall s \in \mathcal{S},\, t \in \mathcal{T}}
  \addConstraint{u_{st}}{\in \{0, 1\}}{\hspace{6em}\forall s \in \mathcal{S},\, t \in \mathcal{T}} \label{u_st_binary}
  \addConstraint{\sum_{w \in \mathcal{N}(v)}p_{vwt}}{=\sum_{s \in \mathcal{S}_v}p_{st} - \sum_{b \in \mathcal{B}_v} p_{bt}}{\hspace{6em}\forall v \in \mathcal{V},\, t \in \mathcal{T}}
  \addConstraint{\sum_{w \in \mathcal{N}(v)}q_{vwt}}{=\sum_{s \in \mathcal{S}_v}q_{st} - \sum_{b \in \mathcal{B}_v} q_{bt}}{\hspace{6em}\forall v \in \mathcal{V},\, t \in \mathcal{T}}
\end{maxi!}

For the purposes of introducing the relaxations, the specification above is sufficient. The full specification, including additional elements such as minimum uptime constraints and penalization of bus imbalances, is detailed in Appendix \ref{app:full_formulation}. This base model, containing all market and node-level constraints, is shared across all relaxations. Each relaxation then augments the model with its own ACOPF constraints, defining how voltages, flows, and reference angles are represented to capture the physical feasibility of the network.

Before introducing the convex relaxations, we recall the standard non-convex line-flow-based ACOPF formulation. Let $Y = G + jB \in \mathbb{C}^{|\mathcal V|\times|\mathcal V|}$ denote the bus admittance matrix of the network, and let 
\begin{align}
    V_{vt} = |V_{vt}|e^{j\theta_{vt}}, \label{V_vt_polar}
\end{align} 
be the complex voltage phasor at bus $v \in \mathcal{V}$ and time $t \in \mathcal{T}$ with phase angle $\theta_{vt} \in [0, 2\pi)$.

For each time period $t \in \mathcal{T}$ and each transmission line $(v,w) \in \mathcal{E}$, let
\begin{align}
S_{vwt} = p_{vwt} + j\, q_{vwt}, \label{S_vwt}
\end{align}
denote the complex power flow from bus $v$ to bus $w$, with $p_{vwt}$ and $q_{vwt}$ the corresponding real and reactive components.

For each line $(v,w) \in \mathcal{E}$, the current flowing from bus $v$ to bus $w$ is given by Ohm's law:
\begin{align}
I_{vwt} = Y_{vw}\, (V_{vt} - V_{wt}), \label{I_vwt}
\end{align}
where $Y_{vw} = G_{vw} + jB_{vw}$ is the series admittance of the line, with $G_{vw}$ the conductance and $B_{vw}$ the susceptance.

The complex power flow on the line is then obtained via the standard AC relationship between voltage and current:
\begin{align}
\boxed{
S_{vwt} = V_{vt}\, \overline{I_{vwt}} = V_{vt}\, \overline{Y_{vw}\, (V_{vt} - V_{wt})},
} \label{power_flow_equation}
\end{align}
where $\overline{(\cdot)}$ denotes complex conjugation.
This nonlinear, non-convex constraint is the fundamental equation governing power flow in AC networks.

Similar to the bounds on generation and consumption for sellers and buyers, the network imposes operational limits on voltages and line flows. For each bus $v \in \mathcal{V}$ and time period $t \in \mathcal{T}$, the voltage magnitude is constrained by
\begin{align}
\underline{V}_{v} \le |V_{vt}| \le \overline{V}_{v}, \label{V_vt_bounds}
\end{align}
where $\underline{V}_{v}$ and $\overline{V}_{v}$ denote the minimum and maximum allowable voltages at bus $v$.  

Each transmission line $(v,w) \in \mathcal{E}$ also has a thermal (current) limit, denoted by $\overline{I}_{vw}$, restricting the magnitude of the complex current:
\begin{align}
|I_{vwt}| \le \overline{I}_{vw} \quad \forall t \in \mathcal{T}. \label{I_vwt_bounds}
\end{align}
Finally, one bus $r^* \in \mathcal{V}$ is chosen as the reference bus, with 
\begin{align}
    V_{r^*t}=1, \label{reference_bus}
\end{align}
for all $t \in \mathcal{T}$, fixing both the voltage magnitude to 1 p.u. and the phase angle to zero.

This collection of constraints \eqref{power_flow_equation}--\eqref{reference_bus} defines the classical ACOPF feasible region. The convex formulations we present in the next section replace the nonlinear voltage-flow relations with tractable convex surrogates that either approximate or contain this region.

\section{Convex Relaxations of ACOPF}
\subsection{DC Approximation}
\label{sec:dc}
Although not a relaxation, the DC optimal power flow (DCOPF) is a widely used \emph{linear approximation} of the ACOPF model.
It is often presented alongside convex relaxations as a computational baseline, offering fast solve times but weaker solution guarantees.
The DCOPF model assumes that voltage magnitudes remain close to their nominal values and that phase angle differences across lines are small.

To derive the approximation, we first rewrite \eqref{power_flow_equation} using \eqref{V_vt_polar}. This yields
\[
S_{vwt} = |V_{vt}| e^{j\theta_{vt}} (G_{vw} - j B_{vw}) (|V_{vt}| e^{-j\theta_{vt}} - |V_{wt}| e^{-j\theta_{wt}}).
\]

Separating real and imaginary parts and using $e^{j(\theta_{vt}-\theta_{wt})} = \cos(\theta_{vt}-\theta_{wt}) + j \sin(\theta_{vt}-\theta_{wt})$, we obtain the AC power flow equations in polar form:
\begin{align}
p_{vwt} &= |V_{vt}|^2 G_{vw} - |V_{vt}||V_{wt}| \Big( G_{vw} \cos(\theta_{vt}-\theta_{wt}) + B_{vw} \sin(\theta_{vt}-\theta_{wt}) \Big), \label{p_vwt_polar} \\
q_{vwt} &= -|V_{vt}|^2 B_{vw} - |V_{vt}||V_{wt}| \Big( G_{vw} \sin(\theta_{vt}-\theta_{wt}) - B_{vw} \cos(\theta_{vt}-\theta_{wt}) \Big). \label{q_vwt_polar}
\end{align}
These equations express active and reactive power flows in terms of voltage magnitudes and phase angles.

The DC approximation is obtained by making the following assumptions:

\begin{enumerate}[label=(\roman*)]
  \item Flat voltage magnitudes: $|V_{vt}| \approx 1$ p.u.
  \item Small angle differences: $\theta_{vt} - \theta_{wt} \ll 1$, allowing $\sin(\theta_{vt} - \theta_{wt}) \approx \theta_{vt} - \theta_{wt}$ and $\cos(\theta_{vt} - \theta_{wt}) \approx 1$.
  \item Neglect line resistances: $G_{vw} \approx 0$, assuming lines are primarily inductive.
\end{enumerate}

Under these assumptions, the active power flow \eqref{p_vwt_polar} linearizes to
\begin{align}
p_{vwt} \approx B_{vw} (\theta_{vt} - \theta_{wt}), \label{p_vwt_dc}
\end{align}
while the reactive flow \eqref{q_vwt_polar} disappears.

\subsubsection*{Power flow constraints} The final set of power flow constraints corresponding to \eqref{power_flow_equation}--\eqref{reference_bus} become
\begin{align}
    p_{vwt} &= B_{vw} (\theta_{vt} - \theta_{wt}) &&\forall v \in \mathcal{V},\, w \in \mathcal{N}(v),\, t\in\mathcal{T}, \\
    |p_{vwt}| &\le \overline{I}_{vw} &&\forall v \in \mathcal{V},\, w \in \mathcal{N}(v),\, t \in \mathcal{T}, \label{I_vwt_bounds_DCOPF} \\
    \theta_{r^*t} &= 0 &&\forall t \in \mathcal{T}.
\end{align}
The thermal limit \eqref{I_vwt_bounds} on current magnitude translates to a limit on active power flow under the DC assumptions: with $|V_{vt}| \approx 1$ and $|S_{vwt}| \approx |p_{vwt}|$ (since $q_{vwt} = 0$), we have $|I_{vwt}| \approx |p_{vwt}|$, yielding \eqref{I_vwt_bounds_DCOPF}.

\subsubsection*{Voltage reconstruction} Since voltage magnitudes are fixed at 1 p.u. in the DC approximation, we construct the complex voltages from the solved phase angles via \eqref{V_vt_polar}:
\begin{align}
    V_{vt} := \exp(j\theta_{vt}). \label{V_vt_dc}
\end{align}

\subsection{Shor's SDP Relaxation}
\label{sec:shor}
Shor's semidefinite programming (SDP) relaxation is a convex relaxation of the ACOPF that replaces the non-convex voltage product terms with a lifted matrix representation. 
This transformation isolates the non-convexity into a single rank-one constraint, which we then relax to obtain a tractable semidefinite program.

To derive this relaxation, we first rewrite the AC power flow equation \eqref{power_flow_equation} in a convenient form:
\begin{align}
    S_{vwt} = \overline{Y_{vw}} \,(V_{vt} \overline{V_{vt}} - V_{vt} \overline{V_{wt}}). \label{S_vwt_lifted_form}
\end{align}
Here, it is evident that the source of non-convexity is the interaction term $V_{vt} \overline{V_{wt}}$. To address this, we introduce the complex lifted voltage matrix
\begin{align}
\mathbf{W}_t = \mathbf{V}_t \mathbf{V}_t^\mathrm{H} \in \mathbb{C}^{|\mathcal{V}| \times |\mathcal{V}|}, \label{W_rank1_complex_def}
\end{align}
where $\mathbf{V}_t = \bigl(V_{vt}\bigr)_{v \in \mathcal{V}} \in \mathbb{C}^{|\mathcal{V}|}$ is the vector of bus voltage phasors at time $t \in \mathcal{T}$.

For convenience, we write $\mathbf{W}_t = (W_{vwt})_{v,w \in \mathcal{V}}$. Then \eqref{S_vwt_lifted_form} can be expressed as
\begin{align}
    S_{vwt} = \overline{Y_{vw}} \,(W_{vvt} - W_{vwt}), \label{S_vwt_linear_in_W}
\end{align}
which is linear in $W_{vvt}$ and $W_{vwt}$. Thus, the non-convexity in the original power flow is now isolated in the rank-one constraint \eqref{W_rank1_complex_def}. The relaxation is achieved by replacing \eqref{W_rank1_complex_def} with a positive semidefiniteness constraint
\begin{align}
    \mathbf{W}_t &\succeq 0, \label{W_t_psd} \\
    \mathbf{W}_t &= \mathbf{W}_t^\mathrm{H},
\end{align}
where the Hermitian constraint is included explicitly to ensure that $\mathbf{W}_t$ itself is PSD. Some SDP solvers may otherwise interpret the PSD constraint on a complex matrix $\mathbf{W}_t$ as requiring only that its Hermitian part $\frac{1}{2}(\mathbf{W}_t + \mathbf{W}_t^\mathrm{H})$ be PSD.

For the thermal limit \eqref{I_vwt_bounds}, we reformulate the constraint to obtain a tractable formulation. Our goal is to obtain a constraint that can be expressed as an SOCP without introducing additional variables.  

Rewriting \eqref{I_vwt_bounds} using \eqref{S_vwt} and \eqref{power_flow_equation} gives
\begin{align}
    |S_{vwt} V_{vt}^{-1}| \le \overline{I}_{vw}. \label{thermal_intermediate}
\end{align}

One might consider introducing the lifted matrix $\mathbf{W}_t$ by squaring both sides; however, this leads to a left-hand side that is bilinear in $S_{vwt}$, which cannot be represented as an SOCP constraint. To obtain an SOCP-representable constraint, we leverage the lower bound on $|V_{vt}|$. Since $|V_{vt}| \geq \underline{V}_v$, the constraint \eqref{thermal_intermediate} is satisfied if $|S_{vwt}| \leq \overline{I}_{vw} \underline{V}_v$ holds.  

By rewriting $S_{vwt}$ as a 2-dimensional vector, we obtain the following conservative but SOCP-representable constraint:
\begin{align}
    \left\lVert 
    \begin{bmatrix}
        p_{vwt} \\
        q_{vwt}
    \end{bmatrix}
    \right\rVert_2 
    \le \overline{I}_{vw} \, \underline{V}_{v}. \label{SOCP_thermal}
\end{align}

\subsubsection*{Power flow constraints} The full set of relaxed power flow constraints are
\begin{align}
    \mathbf{W}_t &\succeq 0 && \forall t \in \mathcal{T},\\
    \mathbf{W}_t &= \mathbf{W}_t^\mathrm{H} && \forall t \in \mathcal{T},\\
    p_{vwt} &= \Re(\overline{Y_{vw}}(W_{vvt} - W_{vwt})) &&\forall v \in \mathcal{V},\, w \in \mathcal{N}(v), \, t \in \mathcal{T}, \\
    q_{vwt} &= \Im(\overline{Y_{vw}}(W_{vvt} - W_{vwt})) &&\forall v \in \mathcal{V},\, w \in \mathcal{N}(v), \, t \in \mathcal{T}, \\
        \left\lVert 
    \begin{bmatrix}
        p_{vwt} \\
        q_{vwt}
    \end{bmatrix}
    \right\rVert_2 
    &\le \overline{I}_{vw} \, \underline{V}_{v} &&\forall v \in \mathcal{V},\, w \in \mathcal{N}(v), \, t \in \mathcal{T}, \\
    \underline{V}_{v}^2 &\le \Re(W_{vvt}) \le \overline{V}_{v}^2 &&\forall v \in \mathcal{V},\, t \in \mathcal{T}, \label{shor_complex_V_bounds} \\
    W_{r^*r^*t} &= 1 &&\forall t \in \mathcal{T}, \label{shor_complex_ref0} \\ 
    \Im(W_{r^*vt}) &= 0 &&\forall v \in \mathcal{V}, \, \forall t \in \mathcal{T}, \\
    \Im(W_{vr^*t}) &= 0 &&\forall v \in \mathcal{V}, \, \forall t \in \mathcal{T}. \label{shor_complex_ref1}
\end{align}
The constraint \eqref{shor_complex_V_bounds} corresponds to the original voltage bounds \eqref{V_vt_bounds}. Note that diagonal elements $W_{vvt} = |V_{vt}|^2$ are real for Hermitian matrices, but we retain $\Re(\cdot)$ for numerical robustness. The constraint \eqref{shor_complex_ref0} enforces unit magnitude at the reference bus. The imaginary part constraints in \eqref{shor_complex_ref1} enforce that $V_{r^*t}$ is real (zero phase angle), since $W_{r^*vt} = V_{r^*t}\overline{V_{vt}}$ and requiring $\Im(W_{r^*vt}) = 0$ for all $v \in \mathcal{V}$ implies $V_{r^*t} \in \mathbb{R}$.

\subsubsection*{Voltage reconstruction} In general, the optimal solution $\mathbf{W}_t$ obtained from the relaxation will not be rank-one. To approximate a voltage vector $\mathbf{V}_t$, we project $\mathbf{W}_t$ back to a rank-one representation.  

First, we symmetrize the solution matrix to remove numerical asymmetries:
\[
\widetilde{\mathbf{W}}_t := \tfrac{1}{2}\bigl(\mathbf{W}_t + \mathbf{W}_t^\mathrm{H}\bigr).
\]

Next, we compute the eigenvalue decomposition
\[
\widetilde{\mathbf{W}}_t = \sum_{k=1}^{|\mathcal{V}|} \lambda_k \, u_k u_k^\mathrm{H},
\]
where $\lambda_1 \geq \lambda_2 \geq \dots \geq 0$ are the eigenvalues and $u_k$ the corresponding unit eigenvectors.  

We approximate $\widetilde{\mathbf{W}}_t$ by its dominant rank-one component:
\[
\widetilde{\mathbf{W}}_t \approx \lambda_1 u_1 u_1^\mathrm{H}.
\]

From this representation, we define the approximate bus voltage vector
\[
\widehat{\mathbf{V}}_t := \sqrt{\lambda_1}\, u_1.
\]

Finally, to enforce the reference bus phase constraint, we rotate the recovered vector so that the angle of the reference bus voltage is zero. We set
\[
\mathbf{V}_t := \widehat{\mathbf{V}}_t \, \exp\!\left(-j \,\arg(\widehat{V}_{r^*t})\right),
\]
so that the reference bus has unit magnitude and zero phase angle.

\subsection{Shor's SDP Relaxation (Real-Valued Version)}
The lifted formulation in \eqref{W_rank1_complex_def} is expressed in terms of a Hermitian complex matrix $\mathbf{W}_t \in \mathbb{C}^{|\mathcal{V}|\times|\mathcal{V}|}$.  
While mathematically natural, this complex-valued representation is inconvenient in practice: most standard semidefinite programming solvers operate natively with \emph{real} symmetric matrices rather than complex Hermitian ones.

To bridge this gap, we reformulate the problem using real variables only. Let $v_1, v_2, \ldots, v_{|\mathcal{V}|}$ denote an arbitrary but fixed enumeration of the buses in $\mathcal{V}$.  
We do so by separating real and imaginary parts of the voltage phasors,  
\begin{align}
V_{vt} = V^d_{vt} + j V^q_{vt},
\end{align}
with $V^d_{vt} = |V_{vt}|\cos(\theta_{vt})$ and $V^q_{vt} = |V_{vt}|\sin(\theta_{vt})$
and stacking them into a $2|\mathcal{V}|$-dimensional real vector  
\begin{align}
\widetilde{\mathbf{V}}_t = 
\begin{bmatrix}
V^d_{v_1t} \\[4pt]
V^q_{v_1t} \\
\vdots \\
V^d_{v_nt} \\
V^q_{v_nt}
\end{bmatrix}
\in \mathbb{R}^{2|\mathcal{V}|}.
\end{align}

The corresponding lifted voltage matrix is then defined as
\begin{align}
\mathbf{W}_t = \widetilde{\mathbf{V}}_t \, \widetilde{\mathbf{V}}_t^\top \in \mathbb{R}^{2|\mathcal{V}| \times 2|\mathcal{V}|}, \label{W_rank1_real_def}
\end{align}
which is symmetric, positive semidefinite, and rank-one. 
The matrix $\mathbf{W}_t$ has a block structure where each pair of nodes $(v,w) \in \mathcal{V} \times \mathcal{V}$ corresponds to a $2 \times 2$ block. For node $v$, the diagonal block captures the voltage products at that node:
\begin{align}
    \mathbf{W}_t = \begin{pmatrix}
    \cdot & \cdot & \cdot & \cdots & \cdot & \cdot & \cdot \\ 
    \cdot & (V^d_{vt})^2 & V^d_{vt}V^q_{vt} & \cdots & V^d_{vt}V^d_{wt} & V^d_{vt}V^q_{wt} &\cdot \\
    \cdot & V^q_{vt}V^d_{vt} & (V^q_{vt})^2 & \cdots & V^q_{vt}V^d_{wt} & V^q_{vt}V^q_{wt} & \cdot \\
    \vdots & \vdots & \vdots & \ddots & \vdots & \vdots & \vdots \\
    \cdot & V^d_{wt}V^d_{vt} & V^d_{wt}V^q_{vt} & \cdots & (V^d_{wt})^2 & V^d_{wt}V^q_{wt} & \cdot \\
    \cdot & V^q_{wt}V^d_{vt} & V^q_{wt}V^q_{vt} & \cdots & V^q_{wt}V^d_{wt} & (V^q_{wt})^2 & \cdot \\
    \cdot & \cdot & \cdot & \cdots & \cdot & \cdot & \cdot \\ 
    \end{pmatrix}.
\end{align}
It is a standard result in linear algebra \cite{horn2013matrix, Lavaei2012ZeroDG, Madani2015} that positive semidefiniteness of the Hermitian lifted matrix in \eqref{W_rank1_complex_def} is equivalent to positive semidefiniteness of the real-valued lifted matrix above.

To express the power flow constraints compactly, we introduce shorthand notation for specific matrix elements and their combinations that correspond to physical quantities:
\begin{align*}
    W_{(2v)(2w)t} &= V^d_{vt}V^d_{wt}, \\
    W_{(2v)(2w + 1)t} &= V^d_{vt}V^q_{wt}, \\
    W_{(2v + 1)(2w)t} &= V^q_{vt}V^d_{wt}, \\
    W_{(2v + 1)(2w + 1)t} &= V^q_{vt}V^q_{wt}, \\
    W^\Delta_{vt} &= W_{(2v)(2v)t} + W_{(2v + 1)(2v + 1)t} = |V_{vt}|^2, \\
    W^{\text{Re}}_{vwt} &= W^\Delta_{vt} - (W_{(2v)(2w)t} + W_{(2v + 1)(2w + 1)t}) = |V_{vt}|^2 - \Re(V_{vt}\overline{V_{wt}}), \\
    W^{\text{Im}}_{vwt} &= W_{(2v)(2w + 1)t} - W_{(2v + 1)(2w)t} = \Im(V_{vt}\overline{V_{wt}}).
\end{align*}

To derive the power flow expressions in rectangular coordinates, we start from \eqref{S_vwt_lifted_form}:
\begin{align*}
    S_{vwt} &= \overline{Y_{vw}}(V_{vt}\overline{V_{vt}} - V_{vt}\overline{V_{wt}}) \\
    &= (G_{vw} - jB_{vw})(|V_{vt}|^2 - V_{vt}\overline{V_{wt}}).
\end{align*}
Separating real and imaginary parts and using $\Re(V_{vt}\overline{V_{wt}}) = V^d_{vt}V^d_{wt} + V^q_{vt}V^q_{wt}$ and $\Im(V_{vt}\overline{V_{wt}}) = V^d_{vt}V^q_{wt} - V^q_{vt}V^d_{wt}$ for $V_{vt} = V^d_{vt} + jV^q_{vt}$:
\begin{align}
    p_{vwt} &= G_{vw}|V_{vt}|^2 - G_{vw}\Re(V_{vt}\overline{V_{wt}}) + B_{vw}\Im(V_{vt}\overline{V_{wt}}), \label{p_vwt_intermediate}\\
    q_{vwt} &= -B_{vw}|V_{vt}|^2 + B_{vw}\Re(V_{vt}\overline{V_{wt}}) + G_{vw}\Im(V_{vt}\overline{V_{wt}}). \label{q_vwt_intermediate}
\end{align}
Using the shorthand notation introduced above, we obtain:
\begin{align}
p_{vwt} &= G_{vw}W^{\text{Re}}_{vwt} + B_{vw}W^{\text{Im}}_{vwt}, \label{p_vwt_rectangular} \\
q_{vwt} &= -B_{vw}W^{\text{Re}}_{vwt} + G_{vw}W^{\text{Im}}_{vwt}. \label{q_vwt_rectangular}
\end{align}

\subsubsection*{Power flow constraints} Using the notation we introduced above, the full set of constraints become:
\begin{align}
    \mathbf{W}_t &\succeq 0 && \forall t \in \mathcal{T},\\
    \mathbf{W}_t &= \mathbf{W}_t^\top && \forall t \in \mathcal{T},\\
    p_{vwt} &= G_{vw}W^{\text{Re}}_{vwt} + B_{vw}W^{\text{Im}}_{vwt} &&\forall v \in \mathcal{V},\, w \in \mathcal{N}(v), \, t \in \mathcal{T}, \label{shor_real_p_vwt}\\
    q_{vwt} &= -B_{vw}W^{\text{Re}}_{vwt} + G_{vw}W^{\text{Im}}_{vwt} &&\forall v \in \mathcal{V},\, w \in \mathcal{N}(v), \, t \in \mathcal{T}, \label{shor_real_q_vwt} \\
        \left\lVert 
    \begin{bmatrix}
        p_{vwt} \\
        q_{vwt}
    \end{bmatrix}
    \right\rVert_2 
    &\le \overline{I}_{vw} \, \underline{V}_{v} &&\forall v \in \mathcal{V},\, w \in \mathcal{N}(v), \, t \in \mathcal{T}, \\
    \underline{V}_{v}^2 &\le W^{\Delta}_{vt} \le \overline{V}_{v}^2 &&\forall v \in \mathcal{V},\, t \in \mathcal{T}, \label{shor_real_V_bounds} \\
    W_{(2r^*)(2r^*)t} &= 1 &&\forall t \in \mathcal{T}, \label{shor_real_ref0a} \\ 
    W_{(2r^* + 1)(2r^* + 1)t} &= 0 &&\forall t \in \mathcal{T}, \label{shor_real_ref0b} \\ 
    W_{(2r^* + 1)(2v)t} &= W_{(2r^* + 1)(2v + 1)t} = 0 &&\forall v \in \mathcal{V}, \, \forall t \in \mathcal{T}, \\
    W_{(2v)(2r^* + 1)t} &= W_{(2v + 1)(2r^* + 1)t} = 0 &&\forall v \in \mathcal{V}, \, \forall t \in \mathcal{T}. \label{shor_real_ref1}
\end{align}
Similar to the complex-valued case, the constraints \eqref{shor_real_p_vwt} and \eqref{shor_real_q_vwt} are linear in $\mathbf{W}_t$. The constraint \eqref{shor_real_V_bounds} corresponds to the original voltage bounds \eqref{V_vt_bounds}. Constraints \eqref{shor_real_ref0a}--\eqref{shor_real_ref0b} enforce $V^d_{r^*t} = 1$ and $V^q_{r^*t} = 0$, ensuring the reference bus voltage is real with unit magnitude. Constraint \eqref{shor_real_ref1} enforces that all products involving the imaginary component of the reference bus voltage are zero (i.e., $W_{(2r^*+1)(2v)t} = V^q_{r^*t}V^d_{vt} = 0$, $W_{(2r^*+1)(2v+1)t} = V^q_{r^*t}V^q_{vt} = 0$, and their symmetric counterparts), which is consistent with $V^q_{r^*t} = 0$.

\subsubsection*{Voltage reconstruction} Recovering the voltages is analogous to the complex-valued case. First, we symmetrize the solution matrix to remove numerical asymmetries:
\[
\widetilde{\mathbf{W}}_t := \tfrac{1}{2}(\mathbf{W}_t + \mathbf{W}_t^\top).
\]
We then compute the eigenvalue decomposition and extract the dominant eigenvector $u_1$ corresponding to the largest eigenvalue $\lambda_1$. From the components of $\sqrt{\lambda_1}u_1$, we construct the intermediate voltage vector (where the prime denotes the unrotated reconstruction):
\[
    \widehat{\mathbf{V}}'_t := (\widehat{V}_{(2v)t} + j\widehat{V}_{(2v + 1)t})_{v \in \mathcal{V}} \in \mathbb{C}^{|\mathcal{V}|}.
\]
Finally, we rotate to enforce the reference bus phase constraint:
\[
\mathbf{V}_t := \widehat{\mathbf{V}}'_t \, \exp\!\left(-j \,\arg(\widehat{V}'_{r^*t})\right).
\]

\subsection{Chordal SDP Relaxation}
\label{sec:chordal}

The full positive semidefiniteness constraint in Shor's SDP relaxation above is computationally expensive and often intractable for large-scale instances. The chordal decomposition technique introduced in this section exploits sparsity in the network to replace the full PSD constraint with smaller, tractable constraints while preserving equivalence under suitable conditions. We begin by presenting the required graph-theoretic background before presenting the relaxation:
\begin{itemize}
    \item Given a graph $G = (\mathcal{V}, \mathcal{E})$, a clique is a subset $C \subseteq \mathcal{V}$ such that $(v, w) \in \mathcal{E}$ for all distinct vertices $v, w \in C$. We call a clique maximal, if it is maximal with respect to set inclusion. 
    \item A graph is chordal if every cycle of length at least four contains a chord, i.e., an edge connecting two nonconsecutive vertices of the cycle. Intuitively, chordal graphs resemble triangulations, where long cycles are decomposed into triangles. 
\end{itemize}
The crux of this relaxation is the following theorem:
\vspace{1em}
\begin{theorem}[\citet{fukuda2001}]
    \label{theorem1}
    Let $G = (\mathcal{V}, \mathcal{E})$ be a chordal graph. Let $\mathcal{C} = \{C_1, \dots, C_k\}$ be the collection of maximal cliques of $G$. Let $\mathbf{W} = (W_{vw})_{v,w\in\mathcal{V}} \in \mathbb{C}^{|\mathcal{V}|\times|\mathcal{V}|}$ be a matrix with rows and columns corresponding to vertices of $G$. Then the following statements are equivalent:
    \begin{enumerate}[label=(\roman*), ref=(\roman*)]
        \item $\mathbf{W}_{C_i} = (W_{vw})_{v,w\in C_i}$ is Hermitian positive semidefinite for all $i=1,\dots,k$,
        \item There exists a Hermitian positive semidefinite matrix $\widehat{\mathbf{W}} \in \mathbb{C}^{|\mathcal{V}|\times|\mathcal{V}|}$ such that $W_{vw} = \widehat{W}_{vw}$ for all $v,w \in C_i$ for all $i = 1,\dots,k$.
    \end{enumerate}
\end{theorem}
\vspace{1em}
A stronger formulation of this result and its proof can be found in \citet{fukuda2001, Nakata2003}. This theorem allows us to replace the positive semidefinite (PSD) constraint on the full matrix \eqref{W_t_psd} with PSD constraints on submatrices corresponding to the maximal cliques of the associated graph. However, the network topology $G = (\mathcal{V}, \mathcal{E})$ is not necessarily chordal, which prevents a direct decomposition into clique-based PSD constraints.  

A practical approach to address this is to construct a \emph{chordal extension} of the graph. This can be achieved by performing a symbolic Cholesky decomposition on the adjacency matrix of $G$ (pseudo-code is provided in Appendix~\ref{app:symbolic_cholesky}). Symbolic Cholesky decomposition is a combinatorial procedure that computes the sparsity pattern of the Cholesky factor of a matrix without performing any numerical factorization. In terms of graph theory, eliminating a node in the symbolic Cholesky process corresponds to connecting all its neighbors in the remaining graph, thereby producing a chordal supergraph that preserves the original adjacency structure. More details can be found in \cite{Scott2023}.

Our objective is to minimize \emph{fill-in}, i.e., the number of additional edges introduced in the chordal extension, so that the resulting clique-based PSD decomposition is as sparse as possible. Finding an elimination order that minimizes fill-in is known to be NP-hard \citep{yannakakis1981}. A commonly used heuristic is to first permute the rows and columns of the adjacency matrix using a \emph{minimum degree ordering} and then perform symbolic Cholesky decomposition on the permuted matrix. This approach significantly reduces fill-in in practice while remaining computationally feasible.

Even after restricting our attention to PSD constraints on the clique submatrices, additional consistency conditions must be enforced: the matrices corresponding to different cliques must agree on the entries associated with nodes that are contained in multiple cliques. Imposing all possible overlap constraints can make the problem grow rapidly, particularly when there are a large number of cliques.

To manage this efficiently, we construct the \emph{clique graph} 
\begin{align}
G_\mathcal{C} = (\mathcal{V}_\mathcal{C}, \mathcal{E}_\mathcal{C}), \qquad w(C_i, C_j) = |C_i \cap C_j|
\end{align}
where each node represents a maximal clique, and edges are weighted by the size of the intersection between the corresponding cliques. From this graph, we extract a maximum-weight spanning tree, which we refer to as the \emph{clique tree}:
\begin{align}
T_\mathcal{C} = (\mathcal{V}_\mathcal{C}, \mathcal{E}_T), \qquad \mathcal{E}_T \subseteq \mathcal{E}_\mathcal{C}.
\end{align}
The clique tree satisfies the \emph{running intersection property} (RIP) \citep{lauritzen1996, fukuda2001}, which asserts that for any two cliques \(C_i\) and \(C_j\), their intersection \(C_i \cap C_j\) is contained in every clique along the unique path in the tree connecting \(C_i\) and \(C_j\). For a more detailed discussion of the significance and applications of clique trees, see \citep{lauritzen1988, fukuda2001, koller1999, hara2006, garstka2020cliquegraphbasedmerging}.

This property ensures that overlap constraints need only be imposed along the edges of the clique tree: enforcing equality of the overlapping entries for neighboring cliques automatically guarantees consistency for all other cliques that contain the same nodes, avoiding redundant constraints.

To further improve computational tractability, it is not enough to decompose the PSD constraint into clique submatrices. We must also manage the size of clique intersections. One effective approach is to construct coarser decompositions by merging cliques. A systematic greedy heuristic for this was proposed by \citet{molzahn2013}. The idea is to repeatedly merge pairs of cliques in such a way that the overall size of the semidefinite program, measured by the number of matrix variables plus the number of linking constraints, is reduced as much as possible.

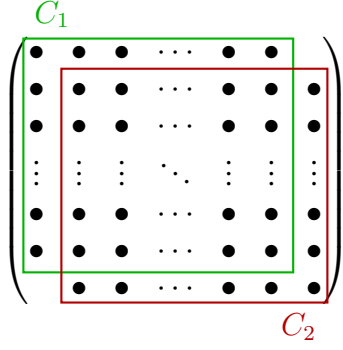
\begin{figure}[H]
\centering
\begin{tikzpicture}[baseline=(current bounding box.center)]
\node (m) {
\(\begin{pmatrix}
\bullet & \bullet & \bullet & \cdots & \bullet & \bullet &  \\
\bullet & \bullet & \bullet & \cdots & \bullet & \bullet & \bullet \\
\bullet & \bullet & \bullet & \cdots & \bullet & \bullet & \bullet \\
\vdots  & \vdots  & \vdots  & \ddots & \vdots  & \vdots  & \vdots  \\
\bullet & \bullet & \bullet & \cdots & \bullet & \bullet & \bullet \\
\bullet & \bullet & \bullet & \cdots & \bullet & \bullet & \bullet \\
 & \bullet & \bullet & \cdots & \bullet & \bullet & \bullet \\
\end{pmatrix}\)
};

\draw[green!70!black,thick] (-2, 1.75) rectangle (1.55, -1.35);
\node[green!70!black,anchor=south west] at (-2, 1.75) {$C_1$};

\draw[red!70!black,thick] (-1.5,1.35) rectangle (2, -1.75);
\node[red!70!black,anchor=north east] at (2, -1.75) {$C_2$};

\end{tikzpicture}
\caption{An example of an unfavorable decomposition. Here two large blocks overlap almost entirely. Decomposing them into separate cliques would introduce a large number of redundant linking constraints. Instead, one should perform the initial chordal decomposition, detect the excessive overlap, and merge the cliques back together.}
\label{clique_merging}
\end{figure}

For bookkeeping, let us define
\[
    \gamma(A) = \frac{|A|(|A|+1)}{2},
\]
which counts the number of complex variables introduced by a clique of size \(|A|\).  
Each clique \(C_i\) contributes \(\gamma(C_i)\) complex variables, while each overlapping pair \((C_i, C_j)\) introduces \(\gamma(C_i \cap C_j)\) linking constraints.  

When two cliques \(C_i\) and \(C_j\) that are adjacent in the clique tree are merged, no extra linking constraints are needed beyond those already implied by the tree. The change in problem size is therefore
\[
    \Delta_{ij} = \gamma(C_i \cup C_j) - \big(\gamma(C_i) + \gamma(C_j) + \gamma(C_i \cap C_j)\big).
\]

The merging procedure then proceeds as follows:  
\begin{enumerate}
    \item Find the edge \((C_i, C_j)\) of the clique tree that minimizes \(\Delta_{ij}\).  
    \item Merge cliques \(C_i\) and \(C_j\).  
    \item Update the clique tree and recompute the corresponding \(\Delta_{ij}\) values for neighboring pairs.  
    \item Repeat until the desired number of cliques remains.  
\end{enumerate}

As a practical guideline, \citet{molzahn2013} suggest merging until the number of cliques is reduced to approximately \(10\%\) of the initial set of maximal cliques. Numerous other heuristics for clique merging exist; for further details, see \citet{garstka2020cliquegraphbasedmerging, Nakata2003}.

\subsubsection*{Power flow constraints} Assume we already have a clique decomposition $\mathcal{C} = \{C_1, \dots, C_k\}$ calculated as described above. Let
\begin{align}
    \mathbf{W}_{C_i, t} = (W_{C_i,vwt})_{v,w \in C_i} \in \mathbb{C}^{|C_i| \times  |C_i|},
\end{align}
denote the clique submatrices for time $t \in \mathcal{T}$.

By construction, every node $v \in \mathcal{V}$ and every edge $(v, w) \in \mathcal{E}$ is contained in at least one maximal clique $C_i$. To write the power flow constraints compactly, we need a mechanism to reference matrix entries that may appear in multiple clique submatrices. We introduce a mapping $\eta$ that assigns to each diagonal entry $(v,v)$ or edge $(v,w)$ a canonical location in the clique decomposition.

Formally, for each node $v \in \mathcal{V}$, we choose some clique $C_{i_v}$ containing $v$ and set $\eta(v,v) = (C_{i_v}, (v,v))$. Similarly, for each edge $(v,w) \in \mathcal{E}$, we choose some clique $C_{i_{vw}}$ containing both $v$ and $w$ and set $\eta(v,w) = (C_{i_{vw}}, (v,w))$. When multiple cliques contain the same node or edge, we arbitrarily select one. The notation $W_{\eta(v,w)t}$ then refers to the entry at position $(v,w)$ in the clique matrix specified by $\eta(v,w)$. Due to the overlap constraints \eqref{overlap_constraint}, any constraint involving this entry automatically ensures consistency across all cliques containing $(v,w)$.

The full set of constraints are:
\begin{align}
    \mathbf{W}_{C_i, t} &\succeq 0 && \forall i \in \{1,\dots,k\},\, t \in \mathcal{T}, \\
    \mathbf{W}_{C_i, t} &= \mathbf{W}_{C_i,t}^\mathrm{H} && \forall i \in \{1,\dots,k\},\, t \in \mathcal{T}, \\
    W_{C_i, vwt} &= W_{C_j,vwt} && \forall (C_i, C_j) \in \mathcal{E}_T,\,\, \{v,w\} \subseteq C_i \cap C_j, \, t \in \mathcal{T}, \label{overlap_constraint} \\
    p_{vwt} &= \Re(\overline{Y_{vw}}(W_{\eta(v,v)t} - W_{\eta(v,w)t})) &&\forall v \in \mathcal{V},\, w \in \mathcal{N}(v), \, t \in \mathcal{T}, \\
    q_{vwt} &= \Im(\overline{Y_{vw}}(W_{\eta(v,v)t} - W_{\eta(v,w)t})) &&\forall v \in \mathcal{V},\, w \in \mathcal{N}(v), \, t \in \mathcal{T}, \\
        \left\lVert 
    \begin{bmatrix}
        p_{vwt} \\
        q_{vwt}
    \end{bmatrix}
    \right\rVert_2 
    &\le \overline{I}_{vw} \, \underline{V}_{v} &&\forall v \in \mathcal{V},\, w \in \mathcal{N}(v), \, t \in \mathcal{T}, \\
    \underline{V}_{v}^2 &\le \Re(W_{\eta(v, v)t}) \le \overline{V}_{v}^2 &&\forall v \in \mathcal{V},\, t \in \mathcal{T}, \label{chordal_V_bounds} \\
    W_{\eta(r^*,r^*)t} &= 1 &&\forall t \in \mathcal{T}, \label{chordal_ref0} \\ 
    \Im(W_{\eta(r^*,v)t}) &= 0 &&\forall v \in \mathcal{V}, \, \forall t \in \mathcal{T}, \\
    \Im(W_{\eta(v,r^*)t}) &= 0 &&\forall v \in \mathcal{V}, \, \forall t \in \mathcal{T}. \label{chordal_ref1}
\end{align}

The clique tree computed earlier is used in \eqref{overlap_constraint} to implement the overlap constraints efficiently. In this section, we formulate all constraints using complex-valued variables. The voltage magnitude and reference bus constraints \eqref{chordal_V_bounds}--\eqref{chordal_ref1} are analogous to those in Shor's relaxation \eqref{shor_complex_V_bounds}--\eqref{shor_complex_ref1}, but applied to the clique-indexed matrix entries via the mapping $\eta$. All of the methods and theorems discussed here also apply to the real-valued case; the only difference is that the indexing becomes considerably more cumbersome. In our implementation, both versions are supported: the CVXPY-based \texttt{ChordalShor} uses complex variables, while the MOSEK implementation operates with real-valued variables.

\subsubsection*{Voltage reconstruction}
After obtaining the positive semidefinite clique matrices $(\mathbf{W}_{C_i, t})_{i=1,\dots,k}$, we assemble them into a global matrix $\mathbf{W}_t \in \mathbb{C}^{|\mathcal{V}|\times|\mathcal{V}|}$.  
For each pair of nodes $v,w \in \mathcal{V}$, the corresponding entry is averaged over all cliques containing both nodes:
\begin{align}
    W_{vwt} 
    = \frac{1}{|\chi(v,w)|} \sum_{i \in \chi(v,w)} W_{C_i,vwt},
\end{align}
where $\chi(v,w) = \{ i \mid v,w \in C_i \}$. Note that by construction of the chordal extension, $\chi(v,w) \neq \emptyset$ for all $(v,w) \in \mathcal{E}$ (every edge is contained in at least one clique). For diagonal entries, $\chi(v,v) \neq \emptyset$ is also guaranteed since every node belongs to at least one clique.
As a consequence of \eqref{overlap_constraint} and the running intersection property, these overlapping entries are guaranteed to agree in exact arithmetic. In practice, however, small numerical discrepancies can occur, and averaging them helps improve numerical stability.

Once we have obtained a partially completed matrix $\mathbf{W}_t$ from the chordal decomposition, the next step is to recover the full voltage matrix. This can be achieved using the proof of Theorem \ref{theorem1} \citep{fukuda2001}, which explicitly shows how to construct the matrix completion.

A central concept in the completion algorithm is the \emph{perfect elimination ordering (PEO)} of the chordal graph $G$, which can be written as
\[
v_{p_1}, v_{p_2}, \dots, v_{p_{|\mathcal{V}|}}.
\]
A PEO is an ordering of the vertices such that, for each $i$, the neighbors of $v_{p_i}$ that appear later in the ordering form a clique:
\begin{align}
    (p_i)_{i=1,\dots,|\mathcal{V}|} \text{ is a PEO } \iff \{v_{p_j} \in \mathcal{N}(v_{p_i}) \mid j > i\} \text{ is a clique for each } i=1,\dots,|\mathcal{V}|.
\end{align}
For chordal graphs, a PEO always exists and can be computed efficiently \citep{golumbic2004}.

The completion algorithm proceeds by iteratively processing the vertices in reverse PEO order. At each step, we impute missing matrix entries using conditional independence relationships implied by the chordal structure. Suppose at iteration $i$ we have already processed the vertices $V = \{v_{p_{|\mathcal{V}|}}, \dots, v_{p_{n-i}}\}$. Define the subsets
\begin{align}
    S &= \{v_{p_{n-i-1}}\}, \\
    U &= V \cap \mathcal{N}(v_{p_{n-i-1}}),\\
    T &= V \setminus U.
\end{align}

Here, $S$ is the single vertex currently being processed, $U$ contains its neighbors among the already-processed vertices, and $T$ contains the remaining processed vertices. By the PEO property, $U$ forms a clique. We partition the matrix $\mathbf{W}_t$ into blocks as follows:
\begin{align}
\mathbf{W}_t = 
\begin{pmatrix}
\ddots & \cdot & \cdot & \cdot \\
\cdot & \mathbf{W}_{SS} & \mathbf{W}_{SU} & \mathbf{W}_{ST} \\
\cdot & \mathbf{W}_{US} & \mathbf{W}_{UU} & \mathbf{W}_{UT} \\
\cdot & \mathbf{W}_{TS} & \mathbf{W}_{TU} & \mathbf{W}_{TT} \\
\end{pmatrix},
\end{align}
where $\mathbf{W}_{AB}$ denotes the submatrix corresponding to rows indexed by $A$ and columns indexed by $B$.

The blocks $\mathbf{W}_{SU}$ and $\mathbf{W}_{US}$ correspond to the connections between the current vertex $v_{p_{n-i-1}} \in S$ and its neighbors $U$. Due to the chordal structure of the graph, these entries are already determined by the PSD constraints on the cliques that contain both $S$ and $U$. 

Since the blocks $\mathbf{W}_{UU}$, $\mathbf{W}_{UT}$, and $\mathbf{W}_{TT}$ involve only previously processed vertices, their entries are already determined. The key observation is that because $U$ forms a clique and separates $S$ from $T$ in the chordal graph, the entries $\mathbf{W}_{ST}$ can be recovered from the Schur complement relationship. Specifically, in a positive semidefinite matrix, the conditional covariance structure allows us to impute the missing entries $\mathbf{W}_{ST}$ as a function of the known blocks. The off-diagonal blocks $\mathbf{W}_{ST}$ and $\mathbf{W}_{TS}$ are updated as
\begin{align}
\mathbf{W}_{ST} &:= \mathbf{W}_{SU} \mathbf{W}_{UU}^\dagger \mathbf{W}_{TU}^\top, \\
\mathbf{W}_{TS} &:= \mathbf{W}_{TU} \mathbf{W}_{UU}^\dagger \mathbf{W}_{SU}^\top,
\end{align}
where $\mathbf{W}_{UU}^\dagger$ denotes the Moore-Penrose pseudoinverse of $\mathbf{W}_{UU}$.

This process continues until all vertices in the PEO have been processed  and guarantees a PSD completion for chordal graphs \citep{grone1984}. After completing the global matrix $\mathbf{W}_t$, voltages are computed as in the standard SDP relaxation. A pseudo-code implementation of the full procedure, including the computation of the PEO, is provided in Appendix \ref{app:psd_completion}.

\subsection{Jabr's Relaxation (SOCP)}
\label{sec:jabr}
Jabr’s relaxation can be viewed as a simplification of Shor's SDP relaxation.  
Instead of lifting the problem into a full matrix variable, it focuses only on edge-based pairwise terms, thereby reducing the semidefinite constraints to second-order cone (SOC) constraints.
We introduce optimization variables $c_{vwt}$ and $s_{vwt}$ to represent the pairwise voltage products and trigonometric terms:
\begin{align}
    c_{vwt} &= |V_{vt}||V_{wt}| \cos(\theta_{vt} - \theta_{wt}), \label{c_vwt_def} \\
    s_{vwt} &= |V_{vt}||V_{wt}| \sin(\theta_{vt} - \theta_{wt}). \label{s_vwt_def}
\end{align}
for $v, w \in \mathcal{V}$ and $t \in \mathcal{T}$. Then, the AC power flow equations in polar form \eqref{p_vwt_polar}, \eqref{q_vwt_polar} become
\begin{align}
p_{vwt} &= c_{vvt} G_{vw} - c_{vwt} G_{vw} - s_{vwt} B_{vw}, \label{p_vwt_jabr} \\
q_{vwt} &= -c_{vvt} B_{vw} - s_{vwt} G_{vw} + c_{vwt} B_{vw}. \label{q_vwt_jabr}
\end{align}
which are fully linear in $c_{vwt}$ and $s_{vwt}$. 

Imposing the constraints
\begin{align}
    c_{vwt} &= c_{wvt}, \\
    s_{vwt} &= -s_{wvt}, \\
    c_{vwt}^2 + s_{vwt}^2 &= c_{vvt}c_{wwt} \label{jabr_nonconvexity}
\end{align}
yields a formulation that is exact for radial networks.
For meshed networks, such as the one considered here, this formulation only represents a relaxation. 
The constraints above do not ensure a consistent set of bus voltage angles whose differences sum to zero (mod $2\pi$) around each cycle in the network \citep{Kocuk2016}. 
Moreover, the constraint \eqref{jabr_nonconvexity} is non-convex. 
To obtain a convex relaxation, we replace it with the rotated SOCP constraint
\begin{align}
    c_{vwt}^2 + s_{vwt}^2 &\le c_{vvt}c_{wwt},
    \label{jabr_rotated_socp}
\end{align}
which can be enforced efficiently within standard SOCP solvers.

Comparing with \eqref{p_vwt_intermediate}--\eqref{q_vwt_intermediate}, we observe that $c_{vwt}$ and $s_{vwt}$ can equivalently be characterized in rectangular coordinates as:
\begin{align}
    c_{vwt} &= \Re(V_{vt}\overline{V_{wt}}) = V^d_{vt}V^d_{wt} + V^q_{vt}V^q_{wt}, \label{c_vwt_rectangular} \\
    s_{vwt} &= \Im(V_{vt}\overline{V_{wt}}) = V^d_{vt}V^q_{wt} - V^q_{vt}V^d_{wt}. \label{s_vwt_rectangular}
\end{align}
This reveals the connection to Shor's SDP relaxation: Jabr's variables $c_{vwt}$ and $s_{vwt}$ correspond to elements of the lifted matrix $\mathbf{W}_t$. Indeed, for each edge $(v,w) \in \mathcal{E}$, the constraint \eqref{jabr_rotated_socp} is equivalent to enforcing positive semidefiniteness of the $2 \times 2$ Hermitian submatrix
\begin{align}
    \begin{bmatrix}
    c_{vvt} & c_{vwt} - js_{vwt} \\
    c_{vwt} + js_{vwt} & c_{wwt}
    \end{bmatrix}
    \succeq 0,
\end{align}
which is precisely the block 
\begin{align}
    \begin{bmatrix} W_{vvt} & W_{vwt} \\ W_{wvt} & W_{wwt} \end{bmatrix} \in \mathbb{C}^{2\times 2}
\end{align}
from Shor's complex matrix $\mathbf{W}_t$. The key difference is that Jabr maintains only edge-based $2 \times 2$ blocks rather than the full $|\mathcal{V}| \times |\mathcal{V}|$ matrix, sacrificing the global coupling enforced by full PSD but gaining computational tractability.

\subsubsection*{Power flow constraints} The complete set of constraints for Jabr’s relaxation is given by:
\begin{align}
    c_{vwt} &= c_{wvt} &&\forall v \in \mathcal{V},\, w \in \mathcal{N}(v), \, t \in \mathcal{T}, \label{jabr_c_vwt_symmetric}\\
    s_{vwt} &= -s_{wvt} &&\forall v \in \mathcal{V},\, w \in \mathcal{N}(v), \, t \in \mathcal{T},\\
            \left\lVert 
    \begin{bmatrix}
        \frac{1}{2}(c_{vvt} - c_{wwt}) \\
        c_{vwt} \\
        s_{vwt} \\
    \end{bmatrix}
    \right\rVert_2 &\le \frac{1}{2}\left(c_{vvt} + c_{wwt}\right) &&\forall v \in \mathcal{V},\, w \in \mathcal{N}(v), \, t \in \mathcal{T}, \label{jabr_rotated_socp_standard}\\
    p_{vwt} &= c_{vvt} G_{vw} - c_{vwt} G_{vw} - s_{vwt} B_{vw} &&\forall v \in \mathcal{V},\, w \in \mathcal{N}(v), \, t \in \mathcal{T}, \label{jabr_p_vwt}\\
    q_{vwt} &= -c_{vvt} B_{vw} - s_{vwt} G_{vw} + c_{vwt} B_{vw} &&\forall v \in \mathcal{V},\, w \in \mathcal{N}(v), \, t \in \mathcal{T}, \label{jabr_q_vwt} \\
        \left\lVert 
    \begin{bmatrix}
        p_{vwt} \\
        q_{vwt}
    \end{bmatrix}
    \right\rVert_2 
    &\le \overline{I}_{vw} \, \underline{V}_{v} &&\forall v \in \mathcal{V},\, w \in \mathcal{N}(v), \, t \in \mathcal{T}, \label{jabr_thermal} \\
    \underline{V}_{v}^2 &\le c_{vvt} \le \overline{V}_{v}^2 &&\forall v \in \mathcal{V},\, t \in \mathcal{T}, \label{jabr_V_bounds}\\
    c_{r^*r^*t} &= 1 &&\forall t \in \mathcal{T}, \label{jabr_ref}
\end{align}
Here, constraint \eqref{jabr_rotated_socp_standard} is the standard form of the rotated SOCP constraint \eqref{jabr_rotated_socp} introduced previously. For the thermal limit, we once again use the SOCP constraint  \eqref{jabr_thermal} which was introduced in \eqref{SOCP_thermal}. 

\subsubsection*{Voltage reconstruction}

After solving the relaxation, the optimal solution provides only the interaction terms $c_{vwt}$ and $s_{vwt}$. To recover voltage magnitudes and phase angles at each bus, we first fix the values for the reference bus:
\[
\theta_{r^*t} := 0, \quad |V_{r^*t}| := 1.
\]

We then select a spanning tree of the network (e.g., via breadth-first search) rooted at the reference bus. Traversing the tree, we propagate voltages from parent nodes to their children: for any edge $(v, w)$ where the voltage of node $v$ is known, the voltage of node $w$ can be computed from \eqref{c_vwt_def} and \eqref{s_vwt_def}.
Solving for $|V_{wt}|$ and $\theta_{wt}$ gives
\begin{align}
    |V_{wt}| &:= \frac{\sqrt{c_{vwt}^2 + s_{vwt}^2}}{|V_{vt}|}, \\
    \theta_{wt} &:= \theta_{vt} - \arctan\!2(s_{vwt}, c_{vwt}),
\end{align}
where $\arctan\!2(y, x)$ is the two-argument arctangent function that returns the angle in the correct quadrant. Repeating this procedure for all edges in the spanning tree allows us to recover voltage magnitudes and phase angles for every bus in the network.

\subsection{Quadratic Convex (QC) Relaxation}
\label{sec:qc}

The QC relaxation builds on Jabr’s SOCP formulation by providing a tighter link between the auxiliary variables and the underlying voltage magnitudes and angle differences. In Jabr’s relaxation, the power flow on line $(v, w)$ is expressed using the auxiliary variables $c_{vwt}$ and $s_{vwt}$ (\eqref{p_vwt_jabr} and \eqref{q_vwt_jabr}). While this approach captures voltage interactions via second-order cone constraints, it does not explicitly account for their nonlinear coupling with the trigonometric functions $\cos(\theta_{vt} - \theta_{wt})$ and $\sin(\theta_{vt} - \theta_{wt})$ as in \eqref{c_vwt_def} and \eqref{s_vwt_def}. Consequently, the relaxation can be relatively loose.

To address this, the QC relaxation introduces additional auxiliary variables that model these nonlinear couplings and replaces the resulting non-convex relationships with \emph{convex envelopes}. Given a function $f$ and a domain $D$, a convex envelope of $f$ over $D$ is a convex set that contains the graph of $f$ over $D$. 

For each \(v \in \mathcal{V}\), $w \in \mathcal{N}(v)$ and time period \(t \in \mathcal{T}\), we define new auxiliary variables $\tilde{\theta}_{vt} \in [0, 2\pi)$, $|\widetilde{V}_{vt}| \in \mathbb{R}_{\ge 0}$ and 
\begin{align}
c'_{vwt} &= \cos(\tilde{\theta}_{vt} - \tilde{\theta}_{wt}), \label{qc_c_vwt_def} \\
s'_{vwt} &= \sin(\tilde{\theta}_{vt} - \tilde{\theta}_{wt}), \label{qc_s_vwt_def} \\
m_{vwt} &= |\widetilde{V}_{vt}| |\widetilde{V}_{wt}|, \label{m_vwt_def}
\end{align}
so that the original auxiliary variables in Jabr's relaxation can be expressed as
\begin{align}
c_{vvt} &= |\widetilde{V}_{vt}|^2, \label{qc_derived_0} \\
c_{vwt} &= m_{vwt} \, c'_{vwt}, \label{qc_derived_1} \\
s_{vwt} &= m_{vwt} \, s'_{vwt}. \label{qc_derived_2}
\end{align}

These expressions introduce four types of non-convex functions that must be replaced with convex envelopes to maintain convexity: the sine function \eqref{qc_s_vwt_def}, the cosine function \eqref{qc_c_vwt_def}, the quadratic term \eqref{qc_derived_0}, and the bilinear terms \eqref{m_vwt_def}, \eqref{qc_derived_1}, and \eqref{qc_derived_2}.

For a nonlinear relationship $z = f(x)$ with $x \in [x_L, x_U]$, or a bilinear relationship $z = xy$ with $x \in [x_L, x_U]$ and $y \in [y_L, y_U]$, we define the convex envelope as a set of linear constraints that outer-approximate the nonlinear relationship.

We denote the convex envelope of a function $f$ over $[x_L,x_U]$ by $\langle x;[x_L,x_U]\rangle^{f}$, and the McCormick envelope of a bilinear term over $[x_L,x_U]\times[y_L,y_U]$ by $\langle x,y;[x_L,x_U]\times[y_L,y_U]\rangle^{\text{mc}}$. Following \citet{coffrin2016} and \citet{hijazi2017}, we define the convex envelopes as
\begin{align}
z \in \langle x; [x_L,x_U]\rangle^{\sin} &\iff
\left\{
\begin{aligned}
    z &\le \cos\Big(\frac{x_M}{2}\Big) \Big(x - \frac{x_M}{2}\Big) + \sin\Big(\frac{x_M}{2}\Big), \\
    z &\ge \cos\Big(\frac{x_M}{2}\Big) \Big(x + \frac{x_M}{2}\Big) - \sin\Big(\frac{x_M}{2}\Big)
\end{aligned}
\right\}, \label{sin_envelope} \\
z \in \langle x; [x_L,x_U]\rangle^{\cos} &\iff
\left\{
\begin{aligned}
    z &\le 1 - \frac{1 - \cos(x_M)}{x_M^2} x^2, \\
    z &\ge \frac{\cos(x_L) - \cos(x_U)}{x_L - x_U} (x - x_L) + \cos(x_L)
\end{aligned}
\right\}, \label{cos_envelope} \\
z \in \langle x; [x_L,x_U]\rangle^{\text{sq}} &\iff
\left\{
\begin{aligned}
    z &\le (x_L+x_U)x - x_L x_U, \\
    z &\ge x^2
\end{aligned}
\right\}, \label{square_envelope} \\
z \in \langle x,y; [x_L,x_U], [y_L,y_U]\rangle^{\text{mc}} &\iff
\left\{
\begin{aligned}
    z &\ge x_L y + y_L x - x_L y_L, \\
    z &\ge x_U y + y_U x - x_U y_U, \\
    z &\le x_L y + y_U x - x_L y_U, \\
    z &\le x_U y + y_L x - x_U y_L
\end{aligned}
\right\}, \label{mccormick_envelope}
\end{align}
where $x_M = \max\{|x_L|, |x_U|\}$.

Geometrically, each envelope defines a convex set that contains all points on the graph of the corresponding nonlinear function over its domain. 
\begin{itemize}
    \item[\eqref{sin_envelope}] For $\sin(x)$, the envelope is a convex region bounded by tangent lines at selected points, which closely approximates the shape of the sine function over the domain.
    \item[\eqref{cos_envelope}] For $\cos(x)$, the envelope consists of an upper bound given by a quadratic function and a lower bound given by a linear function connecting the endpoints, forming a convex set containing the cosine curve.
    \item[\eqref{square_envelope}] For $x^2$, the envelope consists of a linear function forming the upper bound and the quadratic function $x^2$ forming the lower bound.
    \item[\eqref{mccormick_envelope}] For bilinear terms $xy$, the McCormick envelope forms a convex set containing all feasible $(x,y,z)$ triples, providing a tight convex relaxation of the bilinear relationship over the given bounds.
\end{itemize}

The benefit of these convex envelopes relies on having meaningful bounds for the variables; if the bounds are loose or trivial, the QC relaxation does not improve significantly over the original Jabr relaxation.

\subsubsection*{Power flow constraints} Let 
$\Delta\theta_L$, $\Delta\theta_U$, $c'_L$, $c'_U$, $s'_L$, and $s'_U$ be lower and upper bounds for voltage angle differences and their corresponding sine and cosine values, i.e.
\begin{align}
    \theta_{vt} - \theta_{wt} &\in [(\Delta\theta_L)_{vw},(\Delta\theta_U)_{vw}] &\forall v \in \mathcal{V},\, w \in \mathcal{N}(v), \, t \in \mathcal{T}, \\
    \cos(\theta_{vt} - \theta_{wt}) &\in [(c'_L)_{vw},(c'_U)_{vw}] &\forall v \in \mathcal{V},\, w \in \mathcal{N}(v), \, t \in \mathcal{T}, \\
    \sin(\theta_{vt} - \theta_{wt}) &\in [(s'_L)_{vw},(s'_U)_{vw}] &\forall v \in \mathcal{V},\, w \in \mathcal{N}(v), \, t \in \mathcal{T},
\end{align}
for all admissible phase angle values $\theta_{vt}, \theta_{wt}$ for $(v, w) \in \mathcal{E}$.

The constraints for the QC relaxation is given by:
\begin{align}
    c'_{vwt} &\in \langle \tilde{\theta}_{vt} - \tilde{\theta}_{wt}; [(\Delta\theta_L)_{vw},(\Delta\theta_U)_{vw}]\rangle^{\cos} &&\forall v \in \mathcal{V},\, w \in \mathcal{N}(v), \, t \in \mathcal{T}, \label{qc_envelope_constraints0}\\
    s'_{vwt} &\in \langle \tilde{\theta}_{vt} - \tilde{\theta}_{wt}; [(\Delta\theta_L)_{vw},(\Delta\theta_U)_{vw}]\rangle^{\sin} &&\forall v \in \mathcal{V},\, w \in \mathcal{N}(v), \, t \in \mathcal{T},\\
    m_{vwt} &\in \langle |\widetilde{V}_{vt}|,|\widetilde{V}_{wt}|; [\underline{V}_{v},\overline{V}_{v}],[\underline{V}_{w},\overline{V}_{w}]\rangle^{\text{mc}} &&\forall v \in \mathcal{V},\, w \in \mathcal{N}(v), \, t \in \mathcal{T},\\
    c_{vvt} &\in \langle |\widetilde{V}_{vt}|; [\underline{V}_{v},\overline{V}_{v}]\rangle^{\text{sq}} &&\forall v \in \mathcal{V}, \, t \in \mathcal{T},\\
    c_{vwt} \in \langle m_{vwt},c'_{vwt}&; [\underline{V}_{v}\underline{V}_{w},\overline{V}_{v}\overline{V}_{w}], [(c'_L)_{vw},(c'_U)_{vw}]\rangle^{\text{mc}} &&\forall v \in \mathcal{V},\, w \in \mathcal{N}(v), \, t \in \mathcal{T},\\
    s_{vwt} \in \langle m_{vwt},s'_{vwt}&; [\underline{V}_{v}\underline{V}_{w},\overline{V}_{v}\overline{V}_{w}], [(s'_L)_{vw},(s'_U)_{vw}]\rangle^{\text{mc}} &&\forall v \in \mathcal{V},\, w \in \mathcal{N}(v), \, t \in \mathcal{T}, \label{qc_envelope_constraints1}\\
    &\text{Constraints \eqref{jabr_c_vwt_symmetric}--\eqref{jabr_ref}.}
    \label{qc_jabr_copy}
\end{align}

The constraints \eqref{qc_envelope_constraints0}--\eqref{qc_envelope_constraints1} correspond to the auxiliary variable relations \eqref{qc_c_vwt_def}--\eqref{qc_derived_2}.

\subsubsection*{Voltage reconstruction}
After solving the relaxation, bus voltages are recovered using the same spanning tree and BFS approach as in the Jabr relaxation. 
The auxiliary variables $\tilde{\theta}_{vt}$, $|\widetilde{V}_{vt}|$, $c'_{vwt}$, $s'_{vwt}$ and $m_{vwt}$ are introduced solely to model the nonlinear coupling; their individual values are not guaranteed to be physically meaningful.

\subsubsection*{Practical Bounds via Quasi-Monte Carlo}

The effectiveness of the QC relaxation depends critically on having sufficiently tight bounds for voltage angle differences and their associated sine and cosine values:
\[
\Delta\theta_L, \, \Delta\theta_U, \, c'_L, \, c'_U, \, s'_L, \, s'_U.
\]

In many practical datasets, such as the ARPA-E test cases, these bounds are not provided explicitly. 
To address this, we employ a quasi-Monte Carlo (QMC) approach using a low-discrepancy Sobol sequence \citep{sobol1967} to empirically estimate feasible ranges for each line. QMC generates more evenly distributed samples than standard random sampling, providing better coverage of the multidimensional input space with fewer points and improving the reliability of the resulting bounds \citep{morokoff1995, Caflisch1998}.

\paragraph{Global Sampling}
We seek to sample feasible voltage configurations for the power network. Each node \( v \in \mathcal{V} \) is associated with a complex voltage phasor, represented by its magnitude and phase angle. Directly sampling the entire network requires drawing points in a \( 2|\mathcal{V}| \)-dimensional space. However, quasi-Monte Carlo sequences, such as Sobol sequences, tend to lose their low-discrepancy advantages as the dimension increases, with practical efficiency often decreasing beyond roughly 30–40 dimensions \citep{JoeKuo2008, Sobol2011}. Therefore, we partition the network into connected subgraphs $\{(\mathcal{V}_l, \mathcal{E}_l)\}_{l=0}^M$ of approximately 12 nodes and sample voltage configurations on these subnetworks.

For each subgraph $ (\mathcal{V}_l, \mathcal{E}_l) $ with $ |\mathcal{V}_l| \approx 12$ nodes, we generate a Sobol sequence \((\xi_i)_{i=1}^{2^d} \subset [0,1]^{2|\mathcal{V}_l|}\), using a distinct seed for each subgraph. For each sample $\xi_i$, we construct voltage magnitudes and phase angles using the admissible voltage ranges
\begin{align}
|V_{vt,i}| &= \underline{V}_v + \xi_{i, 2k_v} \, (\overline{V}_v - \underline{V}_v), \\
\theta_{vt,i} &= \pi \, (2\xi_{i, 2k_v + 1} - 1),
\end{align}
where the indices \( 2k_v, 2k_v + 1 \in \{1, \dots, 2|\mathcal{V}_l|\} \) correspond to the node \( v \) in subgraph \( \mathcal{V}_l \). The relative angle difference between connected nodes is then defined as
\[
\Delta\theta_{vwt,i} = \arctan\!\big(\tan(\theta_{vt,i} - \theta_{wt,i})\big),
\]
which restricts all sampled differences to the range \((-\pi/2, \pi/2)\).

For each sampled configuration, we evaluate the active and reactive power flows $[p_{vwt,i}, q_{vwt,i}]^\top$ using the polar power flow equations \eqref{p_vwt_polar}, \eqref{q_vwt_polar}. Only samples that satisfy the thermal limit 
\[ 
S_{vwt,i} = \sqrt{p_{vwt,i}^2 + q_{vwt,i}^2} \le \overline{S}_{vw}
\] are retained.

Finally, the empirical bounds for each line \((v, w) \in \mathcal{E}\) are determined from the remaining feasible samples:
\begin{align}
\Delta\theta_L &= \min_{i \,\text{feasible}} \Delta\theta_{vwt,i}, &\quad
\Delta\theta_U &= \max_{i \,\text{feasible}} \Delta\theta_{vwt,i}, \\
s'_L &= \min_{i \,\text{feasible}} \sin(\Delta\theta_{vwt,i}), &\quad
s'_U &= \max_{i \,\text{feasible}} \sin(\Delta\theta_{vwt,i}), \\
c'_L &= \min_{i \,\text{feasible}} \cos(\Delta\theta_{vwt,i}), &\quad
c'_U &= \max_{i \,\text{feasible}} \cos(\Delta\theta_{vwt,i}).
\end{align}

This ensures that the bounds employed in the QC relaxation are derived from physically realizable operating points, producing tighter and more meaningful convex relaxations.

\paragraph{Local Sampling}
While global sampling provides reasonable bounds, it may be overly conservative for the specific operating region of interest. To potentially produce tighter convex envelopes, we use a two-stage local sampling approach: first, solve Jabr's relaxation to obtain an initial solution $\{|V_{vt}^{\text{Jabr}}|, \theta_{vt}^{\text{Jabr}}\}_{v \in \mathcal{V}, t \in \mathcal{T}}$; second, sample locally around this solution using range parameters $\varepsilon_V, \varepsilon_{\Delta\theta} \in (0, 1]$. The key difference from global sampling is that, instead of sampling uniformly over the entire feasible space, we sample in a local neighborhood around the Jabr solution for each node and time period:
\begin{align}
    |V_{vt,i}| &:= V_{vt}^{\text{Jabr}} + \cdot \varepsilon_V \cdot \xi_{i, 2k_v} \cdot \frac{1}{2}(\overline{V}_v - \underline{V}_v), \\
    \theta_{vt, i} &:= \theta_{vt}^{\text{Jabr}} + \varepsilon_{\theta} \cdot \pi \cdot (2\xi_{i, 2k_v+1} -1),
\end{align}

For each sample, we compute the power flows and apply the same feasibility filtering and bound extraction steps as in the global case.

Local sampling can, in principle, yield substantially tighter convex envelopes than either global sampling or Jabr's formulation. However, the effectiveness of the resulting bounds depends on the choice of $\varepsilon_V$ and $\varepsilon_{\theta}$. If these parameters are set too small, the envelopes may become overly restrictive and lead to infeasibility; if set too large, the envelopes become loose and the approach offers little improvement over global sampling.

As a practical heuristic, we set $\varepsilon_V = 0.1$ and $\varepsilon_{\theta} = 0.15$ in our experiments, which balances exploration of the local neighborhood with maintaining feasibility. In some edge cases, even these parameters are too tight and the model either becomes infeasible or fails to improve the current phasor errors compared to the Jabr baseline. In those cases we can simply return the Jabr allocation that we have computed in the local sampling step.

It is important to note that QMC-based sampling does not guarantee coverage of the entire feasible domain. Although QMC sequences generally provide more uniform coverage than independent random samples for the same number of points, any finite set of samples necessarily leaves gaps. This means that the overall approach makes the method an approximation of the true QC relaxation, rather than a mathematically valid relaxation in the strict sense. In practice, empirically derived bounds serve as useful heuristics that often tighten relaxations, but they can fail in pathological or poorly sampled regions.

\section{Implementation}

All relaxations were initially prototyped in CVXPY, which provides a high-level and expressive interface for convex optimization. CVXPY supports complex-valued variables natively, simplifying the implementation of ACOPF relaxations, and allows solver-agnostic modeling, which enabled initial experimentation with different solvers including SCS \citep{scs2016}, Clarabel \citep{goulart2024clarabel}, MOSEK, and Gurobi \citep{gurobi}. Moreover, CVXPY relies on Disciplined Convex Programming (DCP) \citep{diamond2016cvxpy} to ensure that the specified optimization problems are convex and solvable.

While these features make CVXPY ideal for rapid prototyping and experimentation, they come at a cost. In particular, the canonicalization steps required to transform high-level problem descriptions into solver-ready form can become a significant bottleneck, especially for large semidefinite constraints. For networks with more than 100 buses, the canonicalization overhead in CVXPY makes SDP-based relaxations computationally intractable, as the canonicalizer struggles to handle the large number of variables and constraints. Additionally, the layers of abstraction involved in solver interfacing add overhead that can negatively affect runtime.

After validating the correctness of our formulations in CVXPY, all relaxations were reimplemented directly using the MOSEK Fusion API. This direct approach eliminates canonicalization overhead and allows efficient handling of large SDP constraints, highlighting the tradeoff between modeling convenience and computational efficiency. All benchmarking results presented in Section~\ref{sec:results} use the MOSEK Fusion implementations.

While the MOSEK Fusion API offers significant performance advantages, one limitation is that it does not support mixed-integer SDPs. To obtain meaningful solutions, we relax the binary constraint \eqref{u_st_binary} to
\begin{align}
    u_{st} &\in [0, 1], \quad \forall s \in \mathcal{S}, \, t \in \mathcal{T},
\end{align}
solve the relaxed problem, and then simultaneously round all $u_{st}$ variables to the nearest integer before re-optimizing the full problem with these values fixed. In all cases considered in our experiments, the re-optimization remained feasible. More sophisticated methods, such as branch-and-bound \citep{land_doig_1960}, are available for handling mixed-integer programs but are beyond the scope of this comparative study focused on evaluating convex relaxation formulations.

To ensure convergence across different network topologies and relaxation formulations, we include penalty terms in the objective function for bus power imbalances and thermal limit violations. These soft constraints are scaled relative to the welfare objective and allow the solver to find solutions even when strict constraints cannot be simultaneously satisfied, which is particularly important when comparing relaxations with different feasible regions. The complete formulation with penalty coefficients is provided in Appendix~\ref{app:full_formulation}.
\section{Results and Discussion}
\label{sec:results}
\subsection{Dataset and Experimental Design}

Our experiments use the ARPA-E Grid Optimization (GO) Challenge 2 dataset \citep{elbert2024arpa}, which provides generator cost curves, load valuations, and operational capacity limits for market-based welfare optimization. The dataset comprises 20 networks ranging from 617 to 31,777 buses, including 16 publicly available synthetic networks and 4 proprietary industry networks.

We focus on the 617-bus network, which includes 94 generators, 404 loads, and 841 transmission lines, as the basis for our benchmarking. Its moderate size makes it well-suited for systematically evaluating the performance of multiple relaxations while keeping computational demands tractable. We use only the 617-bus network from the dataset, as we are primarily interested in studying how relaxation performance scales with network size rather than comparing across different network topologies. While our formulation supports multi-period optimization, we focus on the single-period scenario provided in the dataset. We use the base case network topology and market data, excluding contingencies, transformer settings, and switchable shunts to focus on the core welfare maximization problem.

All experiments were conducted on an AMD EPYC 9254 24-Core Processor, using 4 CPU cores and 96 GB RAM per task. We used MOSEK version 11.0.28. We report two timing metrics:
\begin{itemize}
    \item \emph{Total runtime}: includes all overhead from network sampling, chordal decomposition (for chordal SDP), QMC bound computation (for QC relaxations), problem construction, and solver execution
    \item \emph{Solver time}: measures only the time from problem submission to MOSEK until solution return
\end{itemize}
This distinction allows us to separate formulation complexity from solver performance.

Our goal is to benchmark various relaxations to evaluate how key performance metrics (welfare, runtime, memory usage, and constraint feasibility) scale with network size. To enable systematic comparison, we evaluate subnetworks of varying sizes sampled from the 617-bus network. Subnetworks are generated using a random connected subgraph sampling procedure: starting from a randomly selected bus, we iteratively add connected neighbors until the desired number of nodes is reached. We constrain each subnetwork to contain at least one generator and one load to ensure viable market conditions. This sampling approach ensures that all subnetworks share the same underlying electrical characteristics (line impedances, voltage bounds, etc.) while varying in size and local topology, allowing us to isolate the effect of network scale on relaxation performance.

To account for variability in the random sampling process, we generate 10 independent subnetworks for each target size. This batch size provides a reasonable balance between statistical reliability and computational cost: it is large enough to smooth out sampling artifacts while remaining computationally feasible given the extensive set of relaxations being compared. For each network size, we solve all relaxations on all 10 subnetworks and then average the resulting welfare, runtime, and memory usage metrics to produce representative performance estimates.

We evaluate subnetworks at 20 different sizes: $N \in \{32, 64, 96, \ldots, 608, 617\}$ (increments of 32 nodes). Note that at the full 617-bus size, all 10 samples correspond to the complete network topology, making the averaging redundant but maintaining consistency with our batched evaluation framework. For each network size, we compare the DC approximation, the QC relaxation with local QMC samplings, Jabr's SOCP relaxation, the real-valued Shor SDP relaxation, and the chordal SDP relaxation.

\subsection{Comparative Evaluation of Convex Relaxations}

\begin{figure}[hbtp]
    \centering
    \includegraphics[width=\textwidth]{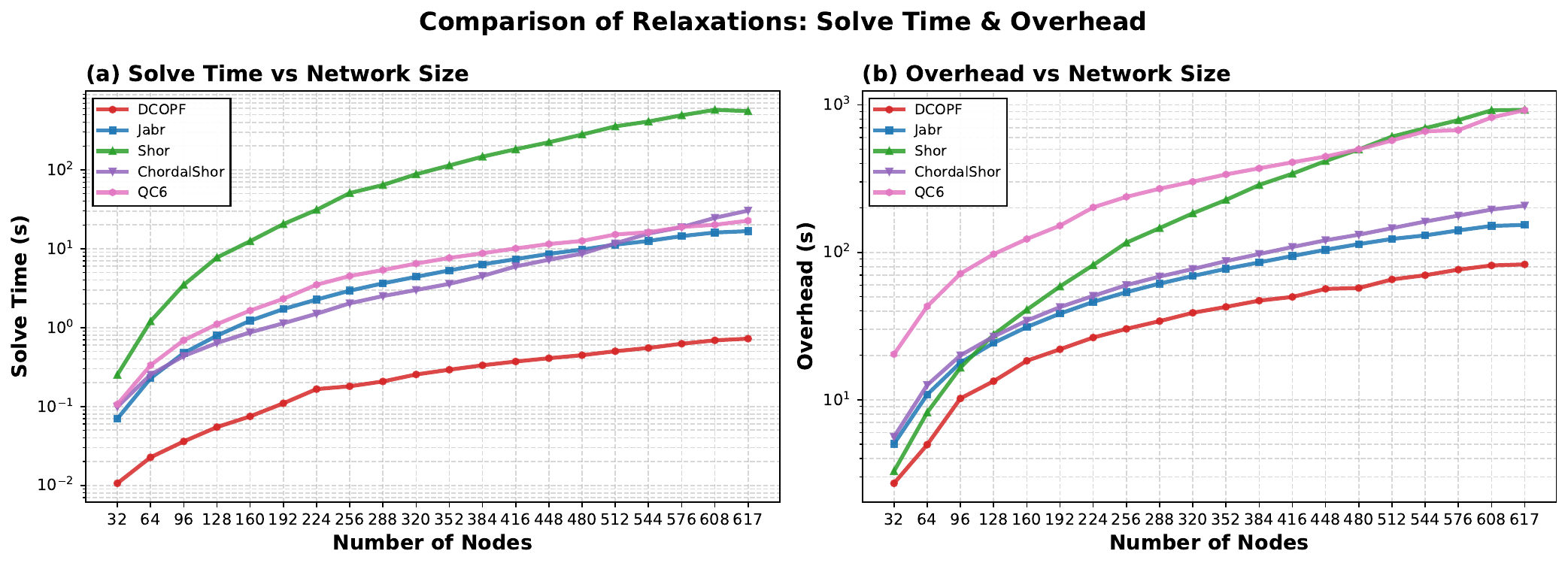}
    \caption{Solve time and overhead (total runtime minus solve time).}
    \label{fig:solve_time_overhead}
\end{figure}

\begin{figure}[H]
    \centering
    \includegraphics[width=\textwidth]{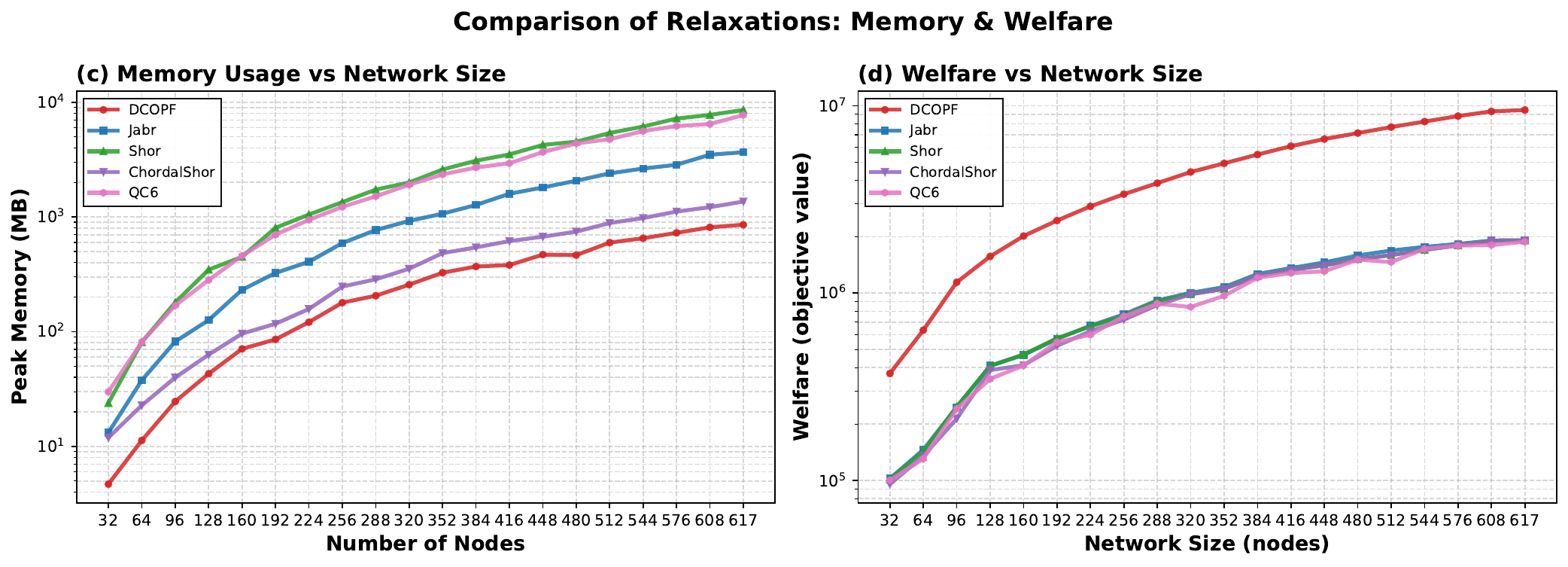}
    \caption{Peak memory usage and economic welfare.}
    \label{fig:memory_welfare}
\end{figure}

\begin{figure}[H]
    \centering
    \includegraphics[width=\textwidth]{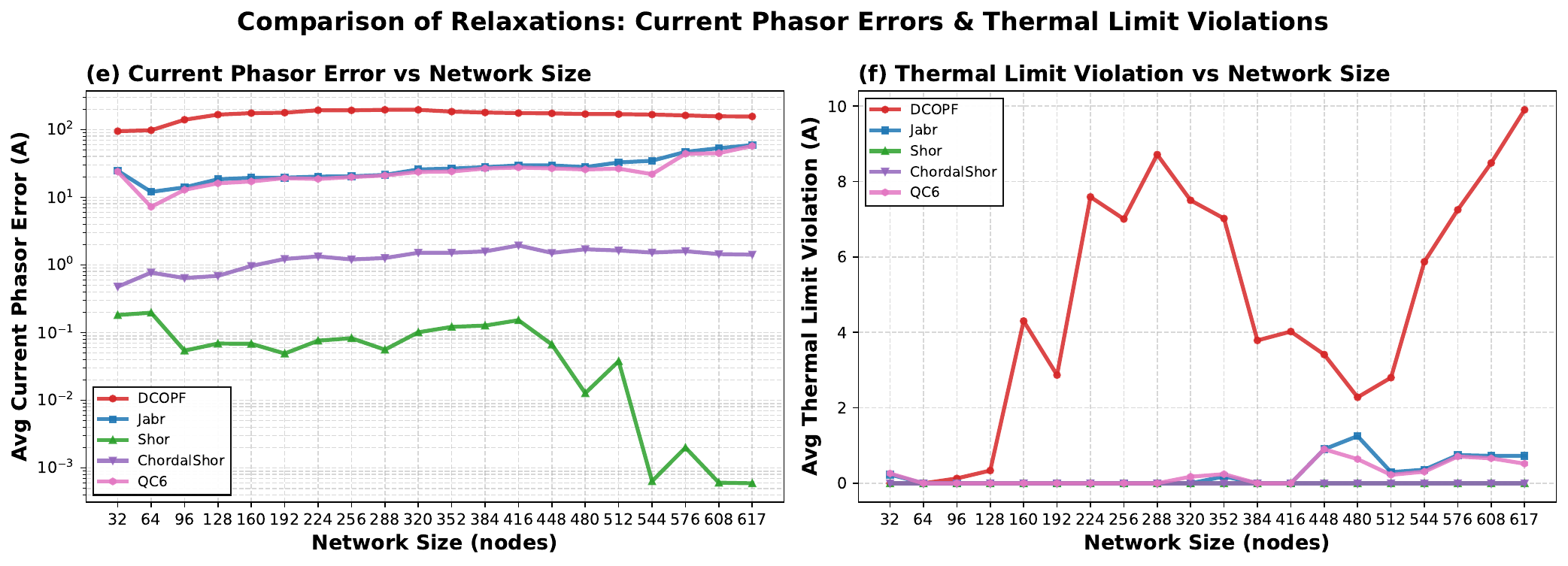}
    \caption{Current phasor error and thermal limit violations (RMS).}
    \label{fig:phasor_violations}
\end{figure}

The practical suitability of each convex relaxation for large-scale market applications depends on both computational efficiency and physical accuracy. We assess performance using six key metrics: \emph{solve time}, \emph{overhead time}, \emph{memory usage}, \emph{welfare}, \emph{current phasor errors}, and \emph{thermal limit violations}. Figures~\ref{fig:solve_time_overhead}–\ref{fig:phasor_violations} summarize the trade-offs between these metrics across all tested methods.

\subsubsection{Solve Time}  
Across all networks, DCOPF is by far the fastest method. Its solve times are typically around one order of magnitude lower than those of Jabr and ChordalShor. Both Jabr and ChordalShor remain highly competitive up to large networks, although ChordalShor tends to scale slightly worse asymptotically. QC6 is marginally slower than both Jabr and ChordalShor because of the extra convex envelope constraints and additional auxiliary variables introduced by the envelope construction. The standard Shor relaxation performs the worst by a wide margin: it is easily an order of magnitude slower than either Jabr or ChordalShor, becoming impractical for large-scale problems.

\subsubsection{Overhead}  
DCOPF also exhibits the smallest preprocessing overhead, as expected. Overall, its total runtime is roughly twice as fast as that of Jabr and ChordalShor. Between the two, ChordalShor incurs slightly higher overhead because it must compute chordal extensions, clique trees, and other graph-structural elements prior to solving. The Shor relaxation, on the other hand, suffers from extreme overhead growth. For small networks its overhead remains moderate—sometimes even smaller than that of Jabr or ChordalShor—but for larger networks it increases dramatically. This is a direct consequence of Shor’s dense semidefinite formulation, which introduces constraints between all pairs of nodes rather than restricting them to existing network edges. Finally, QC6 incurs a lot of extra overhead due to the initial Jabr solution step and the following sampling process.

\subsubsection{Memory Usage}  
Memory consumption follows a clear and intuitive pattern. DCOPF uses the least memory of all methods, followed by ChordalShor, then Jabr. The standard Shor relaxation and QC6 both require substantially more memory. In the case of Shor, this is due to the single large matrix variable that scales quadratically with the number of network nodes. QC6, meanwhile, consumes large amounts of memory because of quasi–Monte Carlo sampling and the storage of the resulting sampled bounds.

\subsubsection{Welfare}  
In terms of market welfare, DCOPF consistently achieves the highest values because it ignores all nonlinear AC feasibility constraints, thereby overestimating total economic surplus. The remaining convex relaxations produce broadly similar welfare levels that are much lower than DCOPF’s. Among these, Jabr and Shor typically perform the best, with nearly identical welfare outcomes. ChordalShor and QC6 yield slightly lower welfare values, indicating that their tighter feasibility enforcement slightly restricts economic efficiency.

\subsubsection{Current Phasor Error}  
The hierarchy of physical accuracy is clearly reflected in the current phasor errors. DCOPF produces the largest errors by a wide margin, demonstrating that its linearized flow model deviates substantially from true AC physics. Jabr performs much better but still exhibits significant deviations. QC6 provides a marginal improvement over Jabr, reducing the average current phasor error slightly. ChordalShor improves accuracy dramatically, with roughly two orders of magnitude lower RMS errors compared to DCOPF. Finally, the standard Shor relaxation achieves exceptionally accurate results, with almost negligible phasor errors, confirming the tightness of the full SDP formulation.

\subsubsection{Thermal Limit Violations}  
Thermal limit violations follow a consistent pattern across all test cases. DCOPF shows persistent and often large violations, especially as network size increases. Jabr and QC6 mitigate most of these violations but not completely. In contrast, ChordalShor, and Shor all satisfy thermal limits strictly, achieving zero or near-zero violations across all tested networks.

\subsection{Comparison of QC Relaxations with Local QMC Sampling Across Sample Counts}

\begin{figure}[H]
    \centering
    \includegraphics[width=\textwidth]{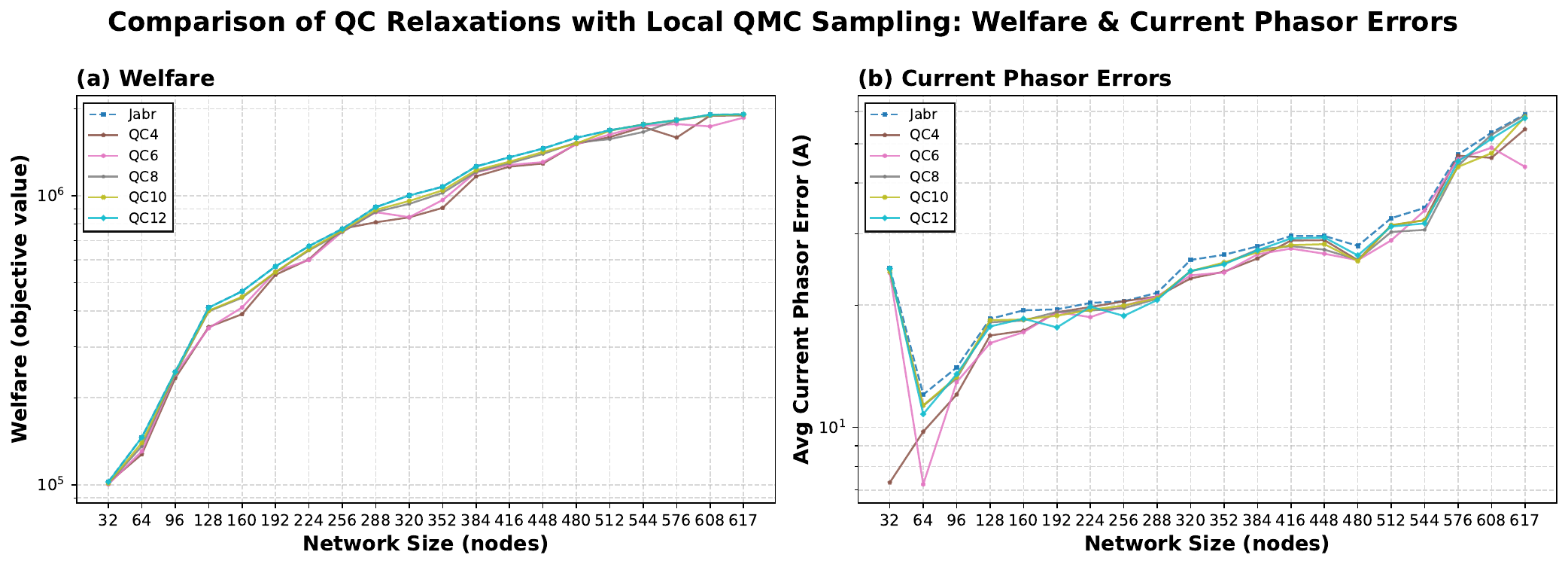}
    \caption{Impact of the number of QMC samples ($2^4$–$2^{12}$) on QC relaxation performance, showing changes in market welfare and current phasor accuracy.}
    \label{fig:qc_degree_plots1}
\end{figure}

\begin{figure}[H]
    \centering
    \includegraphics[width=\textwidth]{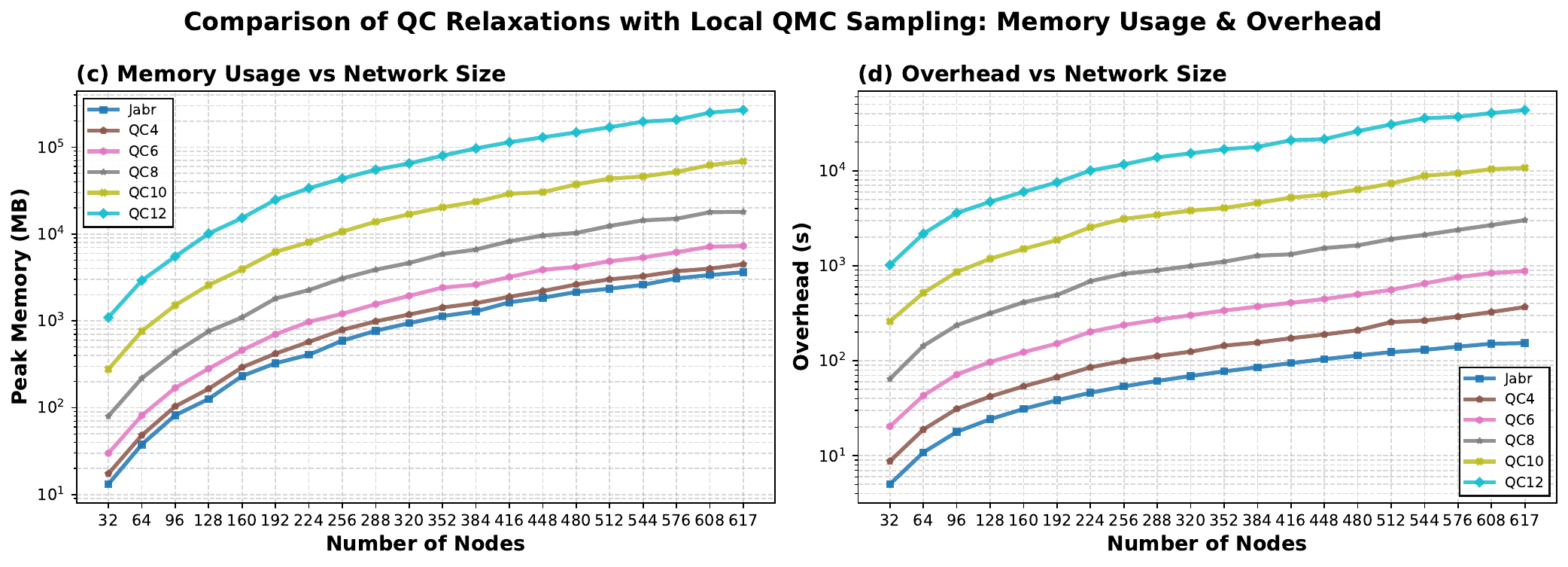}
    \caption{Impact of the number of QMC samples ($2^4$–$2^{12}$) on QC relaxation performance, showing changes in peak memory usage and overhead.}
    \label{fig:qc_degree_plots2}
\end{figure}

All QC variants produce broadly similar welfare outcomes, with all curves remaining slightly below the welfare achieved by the Jabr relaxation. The differences across QC degrees are minor and largely within sampling variability, indicating that increasing the degree does not systematically improve the objective value. QC4 and QC6 sometimes yield marginally lower welfare, but from QC6 onward the results are essentially indistinguishable, suggesting diminishing returns with higher degrees.

The current phasor errors show more variation due to the stochastic nature of the QMC sampling. In some instances, such as QC6 applied to the 64-node subnetworks, the sampled envelopes happen to yield very low current phasor errors. However, this is not a reliable indicator of improved model quality, as infeasibility during sampling can cause fallback to the Jabr solution.

In terms of computational performance, both memory usage and overhead grow rapidly with the degree of the QC relaxation. The observed scaling trends roughly follow the expected $\mathcal O(2^d)$ growth, consistent with the increase in sample sizes. For instance, QC12 requires two to three orders of magnitude more memory and runtime than QC4 for the same network size. On average, QC6 provides the best balance between sampling effort and accuracy improvement, as increasing the degree beyond six offers little to no additional benefit while substantially increasing computational cost.

\subsection{Key Findings}

Our comparative evaluation reveals several clear insights regarding performance, physical accuracy, and practical applicability of the tested convex relaxations:
\begin{itemize}
    \item DCOPF is by far the fastest method, easily outperforming all others by an order of magnitude. Jabr and ChordalShor provide a favorable balance between speed and accuracy, remaining practical even for large networks, while ChordalShor incurs slightly higher overhead due to chordal decomposition. QC6 is somewhat slower due to convex envelope construction, and the standard Shor relaxation is computationally prohibitive for large systems.

    \item DCOPF has the smallest memory footprint, followed by ChordalShor and Jabr. QC6 and Shor require substantially more memory, with QC6’s sampling phase and Shor’s dense SDP matrices driving rapid growth as network size increases.

    \item DCOPF produces the highest nominal welfare by ignoring AC constraints, whereas the physically consistent relaxations (i.e. Jabr, Shor and ChordalShor) yield similar welfare values, with minor reductions reflecting tighter feasibility enforcement.

    \item DCOPF exhibits the largest phasor errors and thermal limit violations. Jabr improves physical fidelity but does not fully eliminate violations. QC6 provides modest further improvement, ChordalShor dramatically reduces errors and achieves full constraint satisfaction, and Shor attains near-perfect accuracy at prohibitive computational cost.

    \item QC relaxations with local QMC sampling are highly sensitive to sampling parameters $(\varepsilon_V, \varepsilon_\theta)$. Proper tuning can improve current phasor accuracy, but infeasible samples revert to the Jabr solution, and stochastic variability leads to inconsistent results. The computational cost of tuning often outweighs the accuracy gains.

    \item For very large networks, DCOPF and Jabr remain the most practical due to speed and scalability. ChordalShor offers the best overall balance between computational tractability and physical fidelity, providing high accuracy, zero violations, and competitive welfare at moderate cost. QC+QMC methods are suitable only for exploratory studies or specialized cases requiring extremely tight current control, while Shor serves primarily as a theoretical benchmark.

\end{itemize}
In summary, ChordalShor is the recommended default relaxation whenever an SDP solver is available, combining robustness, physical realism, and practical efficiency more effectively than any alternative.

\section{Conclusion}

This paper presents a comprehensive comparative evaluation of SDP, SOCP, and QC convex relaxations for market-based AC optimal power flow. Through systematic benchmarking on subnetworks of varying sizes derived from the ARPA-E dataset, we analyze the tradeoffs between solution quality, computational efficiency, and scalability for each relaxation approach.

We made three main contributions. First, we developed a unified framework that specifies and implements the DC approximation as a baseline alongside several widely studied convex relaxations: Shor's SDP relaxation (both complex and real-valued forms), the chordal SDP variant, Jabr's SOCP relaxation, and the QC relaxation. We began by prototyping each model in CVXPY for ease of development, then reimplemented them in the MOSEK Fusion API to reduce modeling overhead and improve computational efficiency, documenting key modeling decisions with particular attention to chordal decomposition techniques and voltage matrix recovery procedures. Second, we addressed the practical challenge that QC relaxations require bounds on voltage angle differences by using quasi-Monte Carlo sampling with Sobol sequences to empirically estimate these bounds when they are unavailable or overly conservative. Third, we presented a comparative study of these relaxations on subnetworks of varying sizes derived from the ARPA-E dataset, assessing not only solution quality but also runtime and memory demands.

Future work could explore higher-order moment relaxations such as the Lasserre hierarchy to assess whether additional tightness justifies the increased computational cost, or apply this comparative framework to other datasets to validate the generality of our findings across different network topologies and market structures.

\section{References}
\bibliographystyle{plainnat}
\bibliography{bibliography}

@book{horn2013matrix,
  title={Matrix Analysis},
  author={Horn, R.A. and Johnson, C.R.},
  isbn={9780521839402},
  lccn={2012012300},
  series={Matrix Analysis},
  url={https://books.google.de/books?id=5I5AYeeh0JUC},
  year={2013},
  publisher={Cambridge University Press}
}

@article{Lavaei2012ZeroDG,
  title={Zero Duality Gap in Optimal Power Flow Problem},
  author={Javad Lavaei and Steven H. Low},
  journal={IEEE Transactions on Power Systems},
  year={2012},
  volume={27},
  pages={92-107},
  url={https://api.semanticscholar.org/CorpusID:2075824}
}

@INPROCEEDINGS{Madani2015,
  author={Madani, Ramtin and Lavaei, Javad and Baldick, Ross},
  booktitle={2015 54th IEEE Conference on Decision and Control (CDC)}, 
  title={Convexification of power flow problem over arbitrary networks}, 
  year={2015},
  volume={},
  number={},
  pages={1-8},
  doi={10.1109/CDC.2015.7402079}}

@ARTICLE{molzahn2017,
  author={Molzahn, Daniel K. and Dörfler, Florian and Sandberg, Henrik and Low, Steven H. and Chakrabarti, Sambuddha and Baldick, Ross and Lavaei, Javad},
  journal={IEEE Transactions on Smart Grid}, 
  title={A Survey of Distributed Optimization and Control Algorithms for Electric Power Systems}, 
  year={2017},
  volume={8},
  number={6},
  pages={2941-2962},
  doi={10.1109/TSG.2017.2720471}}

@article{fukuda2001,
author = {Fukuda, Mituhiro and Kojima, Masakazu and Murota, Kazuo and Nakata, Kazuhide},
title = {Exploiting Sparsity in Semidefinite Programming via Matrix Completion I: General Framework},
journal = {SIAM Journal on Optimization},
volume = {11},
number = {3},
pages = {647-674},
year = {2001},
doi = {10.1137/S1052623400366218},
URL = { https://doi.org/10.1137/S1052623400366218
}
}

@article{waki2006,
author = {Waki, Hayato and Kim, Sunyoung and Kojima, Masakazu and Muramatsu, Masakazu},
title = {Sums of Squares and Semidefinite Program Relaxations for Polynomial Optimization Problems with Structured Sparsity},
journal = {SIAM Journal on Optimization},
volume = {17},
number = {1},
pages = {218-242},
year = {2006},
doi = {10.1137/050623802},
URL = { https://doi.org/10.1137/050623802
}
}

@ARTICLE{jabr2006,
  author={Jabr, R.A.},
  journal={IEEE Transactions on Power Systems}, 
  title={Radial distribution load flow using conic programming}, 
  year={2006},
  volume={21},
  number={3},
  pages={1458-1459},
  doi={10.1109/TPWRS.2006.879234}}

@ARTICLE{coffrin2016,
  author={Coffrin, Carleton and Hijazi, Hassan L. and Van Hentenryck, Pascal},
  journal={IEEE Transactions on Power Systems}, 
  title={The QC Relaxation: A Theoretical and Computational Study on Optimal Power Flow}, 
  year={2016},
  volume={31},
  number={4},
  pages={3008-3018},
  doi={10.1109/TPWRS.2015.2463111}}

@INPROCEEDINGS{hijazi2017,
  author={Coffrin, Carleton and Hijazi, Hassan and Van Hentenryck, Pascal},
  booktitle={2017 IEEE Power \& Energy Society General Meeting}, 
  title={The QC relaxation: A theoretical and computational study on optimal power flow}, 
  year={2017},
  volume={},
  number={},
  pages={1-1},
  doi={10.1109/PESGM.2017.8274004}}

@article{coffrin2021,
title = {The impacts of convex piecewise linear cost formulations on AC optimal power flow},
journal = {Electric Power Systems Research},
volume = {199},
pages = {107191},
year = {2021},
issn = {0378-7796},
doi = {https://doi.org/10.1016/j.epsr.2021.107191},
url = {https://www.sciencedirect.com/science/article/pii/S0378779621001723},
author = {Carleton Coffrin and Bernard Knueven and Jesse Holzer and Marc Vuffray}
}

@BOOK{molzahn2019,
  author={Molzahn, Daniel K. and Hiskens, Ian A.},
  title={A Survey of Relaxations and Approximations of the Power Flow Equations},
  year={2019},
  volume={},
  number={},
  publisher={IEEE},
  pages={},
  doi={10.1561/3100000012}}

@ARTICLE{ghaddar2016,
  author={Ghaddar, Bissan and Marecek, Jakub and Mevissen, Martin},
  journal={IEEE Transactions on Power Systems}, 
  title={Optimal Power Flow as a Polynomial Optimization Problem}, 
  year={2016},
  volume={31},
  number={1},
  pages={539-546},
  doi={10.1109/TPWRS.2015.2390037}}

@article{Kocuk2016,
   title={Strong SOCP Relaxations for the Optimal Power Flow Problem},
   volume={64},
   ISSN={1526-5463},
   url={http://dx.doi.org/10.1287/opre.2016.1489},
   DOI={10.1287/opre.2016.1489},
   number={6},
   journal={Operations Research},
   publisher={Institute for Operations Research and the Management Sciences (INFORMS)},
   author={Kocuk, Burak and Dey, Santanu S. and Sun, X. Andy},
   year={2016},
   month=dec, pages={1177–1196} }

@ARTICLE{matpower2011,
  author={Zimmerman, Ray Daniel and Murillo-Sánchez, Carlos Edmundo and Thomas, Robert John},
  journal={IEEE Transactions on Power Systems}, 
  title={MATPOWER: Steady-State Operations, Planning, and Analysis Tools for Power Systems Research and Education}, 
  year={2011},
  volume={26},
  number={1},
  pages={12-19},
  doi={10.1109/TPWRS.2010.2051168}}

@misc{Lincoln2011,
  author       = {Richard Lincoln},
  title        = {PYPOWER: A Power Flow and Optimal Power Flow Solver in Python},
  year         = {2011},
  howpublished = {\url{https://rwl.github.io/PYPOWER/PYPOWER.pdf}},
  note         = {Accessed: 2025-09-13}
}

@misc{coffrin2018powermodelsjl,
      title={PowerModels.jl: An Open-Source Framework for Exploring Power Flow Formulations}, 
      author={Carleton Coffrin and Russell Bent and Kaarthik Sundar and Yeesian Ng and Miles Lubin},
      year={2018},
      eprint={1711.01728},
      archivePrefix={arXiv},
      primaryClass={math.OC},
      url={https://arxiv.org/abs/1711.01728}, 
}

@INPROCEEDINGS{bitar2012,
  author={Subramanian, A. and Taylor, J. A. and Bitar, E. and Callaway, D. and Poolla, K. and Varaiya, P.},
  booktitle={2012 IEEE 51st IEEE Conference on Decision and Control (CDC)}, 
  title={Optimal power and reserve capacity procurement policies with deferrable loads}, 
  year={2012},
  volume={},
  number={},
  pages={450-456},
  doi={10.1109/CDC.2012.6426102}}

@inbook{kirschen2004,
author    = {Daniel S. Kirschen and Goran Strbac},
title     = {Fundamentals of Power System Economics: Markets for Electrical Energy},
publisher = {John Wiley \& Sons, Ltd},
isbn = {9780470020593},
chapter = {3},
pages = {49-72},
doi = {https://doi.org/10.1002/0470020598.ch3},
url = {https://onlinelibrary.wiley.com/doi/abs/10.1002/0470020598.ch3},
eprint = {https://onlinelibrary.wiley.com/doi/pdf/10.1002/0470020598.ch3},
year = {2004}
}

@inbook{kirschen2004_2,
author    = {Daniel S. Kirschen and Goran Strbac},
title     = {Fundamentals of Power System Economics: Transmission Networks and Electricity Markets},
publisher = {John Wiley \& Sons, Ltd},
isbn = {9780470020593},
chapter = {6},
pages = {141-204},
doi = {https://doi.org/10.1002/0470020598.ch6},
url = {https://onlinelibrary.wiley.com/doi/abs/10.1002/0470020598.ch6},
eprint = {https://onlinelibrary.wiley.com/doi/pdf/10.1002/0470020598.ch6},
year = {2004}
}

@article{moreira2022,
author = {Moreira, Alexandre and Valenzuela, Alan and Heleno, Miguel},
year = {2022},
month = {01},
pages = {1-13},
title = {Solving Market-Based Large-Scale Security-Constrained AC Optimal Power Flows},
volume = {PP},
journal = {IEEE Transactions on Power Systems},
doi = {10.1109/TPWRS.2022.3228211}
}

@ARTICLE{jabr2012,
  author={Jabr, R. A.},
  journal={IEEE Transactions on Power Systems}, 
  title={Exploiting Sparsity in SDP Relaxations of the OPF Problem}, 
  year={2012},
  volume={27},
  number={2},
  pages={1138-1139},
  doi={10.1109/TPWRS.2011.2170772}}

@article{Nakata2003,
  author  = {Kazuhide Nakata and Katsuki Fujisawa and Mituhiro Fukuda and Masakazu Kojima and Kazuo Murota},
  title   = {Exploiting sparsity in semidefinite programming via matrix completion II: implementation and numerical results},
  journal = {Mathematical Programming},
  year    = {2003},
  volume  = {95},
  number  = {2},
  pages   = {303--327},
  doi     = {10.1007/s10107-002-0351-9},
  url     = {https://doi.org/10.1007/s10107-002-0351-9}
}

@article{yannakakis1981,
author = {Yannakakis, Mihalis},
year = {1981},
month = {03},
pages = {},
title = {Computing the Minimum Fill-In is NP-Complete},
volume = {2},
journal = {SIAM Journal on Algebraic and Discrete Methods},
doi = {10.1137/0602010}
}

@Inbook{Scott2023,
author="Scott, Jennifer
and T{\r{u}}ma, Miroslav",
title="Sparse Cholesky Solver: The Symbolic Phase",
bookTitle="Algorithms for Sparse Linear Systems",
year="2023",
publisher="Springer International Publishing",
address="Cham",
pages="53--72",
isbn="978-3-031-25820-6",
doi="10.1007/978-3-031-25820-6_4",
url="https://doi.org/10.1007/978-3-031-25820-6_4"
}

@book{lauritzen1996,
    author = {Lauritzen, Steffen L},
    title = {Graphical Models},
    publisher = {Oxford University Press},
    year = {1996},
    month = {05},
    isbn = {9780198522195},
    doi = {10.1093/oso/9780198522195.001.0001},
    url = {https://doi.org/10.1093/oso/9780198522195.001.0001},
}

@ARTICLE{molzahn2013,
  author={Molzahn, Daniel K. and Holzer, Jesse T. and Lesieutre, Bernard C. and DeMarco, Christopher L.},
  journal={IEEE Transactions on Power Systems}, 
  title={Implementation of a Large-Scale Optimal Power Flow Solver Based on Semidefinite Programming}, 
  year={2013},
  volume={28},
  number={4},
  pages={3987-3998},
  publisher={IEEE},
  doi={10.1109/TPWRS.2013.2258044}
}

@misc{garstka2020cliquegraphbasedmerging,
      title={A clique graph based merging strategy for decomposable SDPs}, 
      author={Michael Garstka and Mark Cannon and Paul Goulart},
      year={2020},
      eprint={1911.05615},
      archivePrefix={arXiv},
      primaryClass={math.OC},
      url={https://arxiv.org/abs/1911.05615}, 
}

@BOOK{andersen2015,
  author={Vandenberghe, Lieven and Andersen, Martin S.},
  title={Chordal Graphs and Semidefinite Optimization},
  year={2015},
  volume={},
  number={},
  pages={},
  journal={Foundations and Trends in Optimization},
  publisher={Now Publishers Inc.},
  doi={10.1561/2400000006}}

@article{lauritzen1988,
    author = {Lauritzen, S. L. and Spiegelhalter, D. J.},
    title = {Local Computations with Probabilities on Graphical Structures and Their Application to Expert Systems},
    journal = {Journal of the Royal Statistical Society: Series B (Methodological)},
    volume = {50},
    number = {2},
    pages = {157-194},
    year = {1988},
    issn = {0035-9246},
    doi = {10.1111/j.2517-6161.1988.tb01721.x},
    url = {https://doi.org/10.1111/j.2517-6161.1988.tb01721.x},
    eprint = {https://academic.oup.com/jrsssb/article-pdf/50/2/157/49097926/jrsssb_50_2_157.pdf},
}

@article{koller1999,
  author       = {Daphne Koller and
                  Uri Lerner and
                  Dragomir Anguelov},
  title        = {A General Algorithm for Approximate Inference and its Application
                  to Hybrid Bayes Nets},
  journal      = {CoRR},
  volume       = {abs/1301.6709},
  year         = {1999},
  url          = {http://arxiv.org/abs/1301.6709},
  eprinttype    = {arXiv},
  eprint       = {1301.6709},
  bibsource    = {dblp computer science bibliography, https://dblp.org}
}

@article{hara2006,
  author       = {Hisayuki Hara and
                  Akimichi Takemura},
  title        = {Boundary cliques, clique trees and perfect sequences of maximal cliques
                  of a chordal graph},
  journal      = {CoRR},
  volume       = {abs/cs/0607055},
  year         = {2006},
  url          = {http://arxiv.org/abs/cs/0607055},
  eprinttype    = {arXiv},
  eprint       = {cs/0607055},
  bibsource    = {dblp computer science bibliography, https://dblp.org}
}

@article{grone1984,
title = {Positive definite completions of partial Hermitian matrices},
journal = {Linear Algebra and its Applications},
volume = {58},
pages = {109-124},
year = {1984},
issn = {0024-3795},
doi = {https://doi.org/10.1016/0024-3795(84)90207-6},
url = {https://www.sciencedirect.com/science/article/pii/0024379584902076},
author = {Robert Grone and Charles R. Johnson and Eduardo M. Sá and Henry Wolkowicz}
}

@incollection{golumbic2004,
  author    = {Martin C. Golumbic},
  title     = {Perfect Elimination Orderings},
  booktitle = {Algorithmic Graph Theory and Perfect Graphs},
  editor    = {Martin C. Golumbic},
  publisher = {Elsevier},
  year      = {2004},
  series    = {Annals of Discrete Mathematics},
  volume    = {57},
  chapter   = {12},
  pages     = {254--266},
  address   = {Amsterdam},
  note      = {Second edition},
}

@article{scs2016,
  author = {O'Donoghue, Brendan and Chu, Eric and Parikh, Neal and Boyd, Stephen},
  title = {Conic Optimization via Operator Splitting and Homogeneous Self-Dual Embedding},
  journal = {Journal of Optimization Theory and Applications},
  year = {2016},
  volume = {169},
  number = {3},
  pages = {1042--1068}
}

@misc{goulart2024clarabel,
  author = {Goulart, Paul and Marmin, Yuwen Chen},
  title = {Clarabel: An interior point conic optimization solver},
  year = {2024},
  url = {https://github.com/oxfordcontrol/Clarabel.rs}
}

@misc{gurobi,
  author = {{Gurobi Optimization, LLC}},
  title = {{Gurobi Optimizer Reference Manual}},
  year = {2023},
  url = {https://www.gurobi.com}
}

@article{diamond2016cvxpy,
  title={CVXPY: A Python-embedded modeling language for convex optimization},
  author={Diamond, Steven and Boyd, Stephen},
  journal={Journal of Machine Learning Research},
  volume={17},
  number={83},
  pages={1--5},
  year={2016}
}

@misc{elbert2024arpa,
  author = {Elbert, S. and Holzer, J. and Veeramany, A. and Coffrin, C. and DeMarco, C. and Duthu, R. and Greene, S. and Kuchar, O. and Lesieutre, B. and Li, H. and Mak, T. and Mittelmann, H. and O'Neill, R. and Overbye, T. and Tbaileh, A. and Van Hentenryck, P. and Wert, J.},
  title = {ARPA-E Grid Optimization (GO) Competition Challenge 2},
  year = {2024},
  publisher = {U.S. Department of Energy},
  url = {https://doi.org/10.25984/2448433}
}

@article{land_doig_1960,
 ISSN = {00129682, 14680262},
 URL = {http://www.jstor.org/stable/1910129},
 author = {A. H. Land and A. G. Doig},
 journal = {Econometrica},
 number = {3},
 pages = {497--520},
 publisher = {[Wiley, Econometric Society]},
 title = {An Automatic Method of Solving Discrete Programming Problems},
 urldate = {2025-10-07},
 volume = {28},
 year = {1960}
}

@article{sobol1967,
author = {I. M. Sobol'},
title  = {Distribution of Points in a Cube
          and the Approximate Evaluation of Integrals (in {R}ussian)},
journal= {Zhurnal Vychislitel'noi Matematiki i Matematicheskoi Fiziki},
year   = 1967,
volume = 7,
pages  = {784--802}
}

@article{morokoff1995,
title = {Quasi-Monte Carlo Integration},
journal = {Journal of Computational Physics},
volume = {122},
number = {2},
pages = {218-230},
year = {1995},
issn = {0021-9991},
doi = {https://doi.org/10.1006/jcph.1995.1209},
url = {https://www.sciencedirect.com/science/article/pii/S0021999185712090},
author = {William J. Morokoff and Russel E. Caflisch}
}

@article{Caflisch1998, title={Monte Carlo and quasi-Monte Carlo methods}, volume={7}, DOI={10.1017/S0962492900002804}, journal={Acta Numerica}, author={Caflisch, Russel E.}, year={1998}, pages={1–49}}

@article{nie2014optimality,
  title={Optimality conditions and finite convergence of Lasserre's hierarchy},
  author={Nie, Jiawang},
  journal={Mathematical Programming},
  volume={146},
  number={1-2},
  pages={97--121},
  year={2014}
}

@article{josz2018lasserre,
  title={Lasserre hierarchy for large scale polynomial optimization in real and complex variables},
  author={Josz, Cédric and Molzahn, Daniel K},
  journal={SIAM Journal on Optimization},
  volume={28},
  number={2},
  pages={1017--1048},
  year={2018}
}

@book{WainwrightJordan2008,
  author    = {Wainwright, Martin J. and Jordan, Michael I.},
  title     = {Graphical models, exponential families, and variational inference},
  publisher = {Now Publishers},
  series    = {Foundations and Trends in Machine Learning},
  volume    = {1},
  year      = {2008}
}

@book{Zhang2006,
  author    = {Zhang, Fanghua},
  title     = {The Schur Complement and Its Applications},
  publisher = {Springer},
  series    = {Numerical Methods and Algorithms},
  volume    = {4},
  year      = {2006}
}

@article{JoeKuo2008,
author = {Joe, Stephen and Kuo, Frances Y.},
title = {Constructing Sobol Sequences with Better Two-Dimensional Projections},
journal = {SIAM Journal on Scientific Computing},
volume = {30},
number = {5},
pages = {2635-2654},
year = {2008},
doi = {10.1137/070709359},
URL = { 
        https://doi.org/10.1137/070709359
},
eprint = { 
        https://doi.org/10.1137/070709359
}
}

@article{Sobol2011,
  author    = {Ilya M. Sobol', Danil Asotsky, Alexander Kreinin},
  title     = {Construction and Comparison of High-Dimensional Sobol' Generators},
  journal   = {Wilmott Journal},
  year      = {2011},
  volume    = {2011},
  pages     = {64--79},
  doi       = {10.1002/wilm.10056},
  url       = {https://onlinelibrary.wiley.com/doi/abs/10.1002/wilm.10056}
}

\appendix

\subsection{Full Formulation}
\label{app:full_formulation}

This appendix contains the complete mathematical formulation of the market-based ACOPF problem we implemented, including penalty terms for constraint violations.

\subsubsection{Decision Variables}
\begin{itemize}
    \item $p_{btl} \ge 0$: power allocated to buyer $b \in \mathcal{B}$ at time $t \in \mathcal{T}$ for block $l \in L_\mathcal{B}$
    \item $p_{stl} \ge 0$: power allocated to seller $s \in \mathcal{S}$ at time $t \in \mathcal{T}$ for block $l \in L_\mathcal{S}$
    \item $p_{bt}$: total active power consumption of buyer $b$ at time $t$
    \item $q_{bt}$: total reactive power consumption of buyer $b$ at time $t$
    \item $p_{st}$: total active power production of seller $s$ at time $t$
    \item $q_{st}$: total reactive power production of seller $s$ at time $t$
    \item $u_{st} \in [0,1]$: unit commitment status of seller $s$ at time $t$ (relaxed from binary)
    \item $\phi_{st} \ge 0$: auxiliary variable for minimum uptime constraints of seller $s$ at time $t$
    \item $p_{vwt}$: active power flow from bus $v$ to bus $w$ at time $t$
    \item $q_{vwt}$: reactive power flow from bus $v$ to bus $w$ at time $t$
    \item $p^\text{imb}_{vt} \ge 0$: active power imbalance at bus $v$ at time $t$ (penalty variable)
    \item $q^\text{imb}_{vt} \ge 0$: reactive power imbalance at bus $v$ at time $t$ (penalty variable)
    \item $I^\text{viol}_{vwt} \ge 0$: thermal limit violation on line $(v,w)$ at time $t$ (penalty variable)
\end{itemize}

Additional variables specific to each relaxation are defined in their respective sections.

\subsubsection{Parameters}

Numerical tolerance and soft constraint parameters:

\begin{itemize}
    \item $\varepsilon_p = 5 \times 10^{-4}$: tolerance for active power flow constraints
    \item $\varepsilon_q = 5 \times 10^{-4}$: tolerance for reactive power flow constraints
    \item $\beta_{p^\text{imb}} = 3 \times 10^{-1}$: scaling for active power balance violations
    \item $\beta_{q^\text{imb}} = 3 \times 10^{-1}$: scaling for reactive power balance violations
    \item $\beta_{I^\text{viol}} = 3 \times 10^{-1}$: scaling for thermal limit violations
\end{itemize}

These parameters enable solver convergence by allowing controlled violations penalized in the objective.

\subsubsection{Objective Function}
The objective maximizes social welfare (buyer valuations minus seller costs and no-load costs) while penalizing constraint violations:

\begin{align}
\max \quad & \sum_{\substack{b \in \mathcal{B} \\ t \in \mathcal{T} \\ l \in L_B}} v_{btl} p_{btl} - \sum_{\substack{s \in \mathcal{S} \\ t \in \mathcal{T} \\ l \in L_S}} c_{stl} p_{stl} - \sum_{\substack{s \in \mathcal{S} \\ t \in \mathcal{T}}} C^\text{NL}_{st} u_{st}  - \alpha_{p^\text{imb}} \sum_{\substack{v \in \mathcal{V} \\ t \in \mathcal{T}}} p^\text{imb}_{vt} - \alpha_{q^\text{imb}} \sum_{\substack{v \in \mathcal{V} \\ t \in \mathcal{T}}} q^\text{imb}_{vt} - \alpha_{I^\text{viol}} \sum_{\substack{(v, w) \in \mathcal{E} \\ t \in \mathcal{T}}} I^\text{viol}_{vwt} \label{obj:full}
\end{align}

where the penalty coefficients are computed as:
\begin{align}
\alpha_\text{welfare} &= \sum_{\substack{b \in \mathcal{B} \\ t \in \mathcal{T} \\ l \in L_B}} |v_{btl}| + \sum_{\substack{s \in \mathcal{S} \\ t \in \mathcal{T} \\ l \in L_S}} |c_{stl}| + \sum_{\substack{s \in \mathcal{S} \\ t \in \mathcal{T}}} |C^\text{NL}_{st}|, \label{welfare_scale} \\
\alpha_{p^\text{imb}} &= \frac{\alpha_\text{welfare} }{|\mathcal{V}| \cdot |\mathcal{T}|}, \label{alpha_p_imb} \\
\alpha_{q^\text{imb}} &= \frac{\alpha_\text{welfare} }{|\mathcal{V}| \cdot |\mathcal{T}|}, \label{alpha_q_imb} \\
\alpha_{I^\text{viol}} &= \frac{\alpha_\text{welfare} }{|\mathcal{V}|^2 \cdot |\mathcal{T}|}. \label{alpha_I_viol}
\end{align}
Here, $\alpha_\text{welfare}$ represents the total scale of the welfare objective, and the penalty coefficients are normalized by the number of corresponding constraint violations to ensure that violations are strongly penalized while maintaining numerical stability and solver convergence.

\subsubsection{Bid and Market Constraints}
\begin{itemize}
    \item \textbf{Buyer constraints:}
\begin{align}
    0 &\le p_{btl} \le \sigma_{btl} &&\forall b \in \mathcal{B},\, t \in \mathcal{T},\, l \in L_\mathcal{B}, \label{cons:buyer_block}\\
    p_{bt} &= \sum_{l \in L_\mathcal{B}} p_{btl} &&\forall b \in \mathcal{B},\, t \in \mathcal{T}, \label{cons:buyer_agg}\\
    \underline{p}_{bt} &\le p_{bt} \le \overline{p}_{bt} &&\forall b \in \mathcal{B},\, t \in \mathcal{T}, \label{cons:buyer_p_bounds}\\
    \underline{q}_{bt} &\le q_{bt} \le \overline{q}_{bt} &&\forall b \in \mathcal{B},\, t \in \mathcal{T}. \label{cons:buyer_q_bounds}
\end{align}

\item \textbf{Seller constraints:}
\begin{align}
    0 &\le p_{stl} \le \sigma_{stl} \cdot u_{st} &&\forall s \in \mathcal{S},\, t \in \mathcal{T},\, l \in L_\mathcal{S}, \label{cons:seller_block}\\
    p_{st} &= \sum_{l \in L_\mathcal{S}} p_{stl} &&\forall s \in \mathcal{S},\, t \in \mathcal{T}, \label{cons:seller_agg}\\
    \underline{p}_{st} \cdot u_{st} &\le p_{st} \le \overline{p}_{st} \cdot u_{st} &&\forall s \in \mathcal{S},\, t \in \mathcal{T}, \label{cons:seller_p_bounds}\\
    \underline{q}_{st} \cdot u_{st} &\le q_{st} \le \overline{q}_{st} \cdot u_{st} &&\forall s \in \mathcal{S},\, t \in \mathcal{T}, \label{cons:seller_q_bounds}\\
    u_{st} &\in [0, 1] \text{ or } u_{st} \in \{0, 1\} &&\forall s \in \mathcal{S},\, t \in \mathcal{T}. \label{cons:unit_commitment}
\end{align}

\item \textbf{Minimum uptime} (multi-period with $|\mathcal{T}| > 1$):
\begin{align}
    \phi_{st} - u_{st} + u_{s,t-1} &\ge 0 &&\forall s \in \mathcal{S},\, t \in \mathcal{T} \setminus \{0\}, \label{cons:phi_def}\\
    \sum_{i=t-T^\text{up}_s}^{t-1} \phi_{si} &\le u_{st} &&\forall s \in \mathcal{S},\, t \in \mathcal{T} : t \ge T^\text{up}_s, \label{cons:min_uptime}
\end{align}
where $T^\text{up}_s$ is the minimum uptime for seller $s$. These constraints ensure that once a generator is started ($u_{st}$ increases from 0 to 1), it must remain online for at least $T^\text{up}_s$ periods.

\end{itemize}

\subsubsection{Network Balance Constraints}

Soft power balance at each bus:
\begin{align}
    \left| \sum_{w \in \mathcal{N}(v)} p_{vwt} - \sum_{s \in \mathcal{S}_v} p_{st} + \sum_{b \in \mathcal{B}_v} p_{bt} \right| &\le p^\text{imb}_{vt} \cdot \beta_{p^\text{imb}} &&\forall v \in \mathcal{V},\, t \in \mathcal{T}, \label{cons:p_balance}\\
    \left| \sum_{w \in \mathcal{N}(v)} q_{vwt} - \sum_{s \in \mathcal{S}_v} q_{st} + \sum_{b \in \mathcal{B}_v} q_{bt} \right| &\le q^\text{imb}_{vt} \cdot \beta_{q^\text{imb}} &&\forall v \in \mathcal{V},\, t \in \mathcal{T}, \label{cons:q_balance}
\end{align}

\subsubsection{Power Flow and Network Constraints}

The power flow constraints, voltage representations, and operational limits depend on the specific relaxation being used:

\begin{itemize}
    \item \textbf{DC Approximation}: Linearized power flow with voltage phase angles only (\ref{sec:dc})
    \item \textbf{Shor's SDP Relaxation}: Lifted matrix formulation with PSD constraints (\ref{sec:shor})
    \item \textbf{Chordal SDP Relaxation}: Clique-based PSD constraints with chordal decomposition (\ref{sec:chordal})
    \item \textbf{Jabr's SOCP Relaxation}: Edge-based auxiliary variables with rotated SOCP constraints (\ref{sec:jabr})
    \item \textbf{QC Relaxation}: Convex envelopes for trigonometric and bilinear terms (\ref{sec:qc})
\end{itemize}

Each relaxation defines voltage variables and power flow equations expressing $p_{vwt}$ and $q_{vwt}$ in terms of these variables. Power flow consistency is enforced with small tolerances $\varepsilon_p$ and $\varepsilon_q$ rather than strict equality. For example, the DC constraint $p_{vwt} = B_{vw}(\theta_{vt} - \theta_{wt})$ becomes $|p_{vwt} - B_{vw}(\theta_{vt} - \theta_{wt})| \le \varepsilon_p$.

\subsubsection{Thermal Limits}
\begin{itemize}
    \item \textbf{DC approximation} (active power only):
\begin{align}
    |p_{vwt}| &\le \overline{I}_{vw} \cdot \left(1 + I^\text{viol}_{vwt} \cdot \beta_{I^\text{viol}}\right) &&\forall v \in \mathcal{V},\, w \in \mathcal{N}(v), \, t \in \mathcal{T}. \label{cons:thermal_dc}
\end{align}

\item\textbf{Other relaxations} (SOCP constraint):
\begin{align}
    \left\lVert 
    \begin{bmatrix}
        p_{vwt} \\
        q_{vwt}
    \end{bmatrix}
    \right\rVert_2 
    &\le \overline{I}_{vw} \underline{V}_v \cdot \left(1 + I^\text{viol}_{vwt} \cdot \beta_{I^\text{viol}}\right) &&\forall v \in \mathcal{V},\, w \in \mathcal{N}(v), \, t \in \mathcal{T}, \label{cons:thermal_socp}
\end{align}

where $\overline{I}_{vw}$ is the line current rating and $\underline{V}_v$ is the minimum voltage magnitude.
 
\end{itemize}

\subsection{Symbolic Cholesky Decomposition}
\label{app:symbolic_cholesky}

\begin{algorithm}[H]
\caption{Symbolic Cholesky Decomposition for Chordal Extension}
\label{alg:symbolic_cholesky}
\begin{algorithmic}[1]
\Require Adjacency matrix $A \in \{0,1\}^{n \times n}$
\Ensure Chordal extension adjacency matrix $A_{\text{chordal}} \in \{0,1\}^{n \times n}$

\State \textcolor{gray}{// Apply minimum degree ordering to reduce fill-in}
\State Compute vertex degrees: $d_v \gets \sum_{w=1}^{n} A[v,w]$ for all $v = 1, \dots, n$
\State Sort vertices by degree: $\pi \gets \text{argsort}(d_1, \dots, d_n)$
\State Permute adjacency matrix: $A \gets A[\pi, \pi]$

\State Initialize chordal extension: $A_{\text{chordal}} \gets A$

\State \textcolor{gray}{// Eliminate vertices in order}
\For{$k = 1$ to $n$}
    \State \textcolor{gray}{// Find uneliminated neighbors of vertex $k$}
    \State $\mathcal{N}_k \gets \{v > k \mid A_{\text{chordal}}[k,v] = 1\}$
    
    \State \textcolor{gray}{// Add fill-ins: connect all neighbors to each other}
    \For{$i \in \mathcal{N}_k$}
        \For{$j \in \mathcal{N}_k$ where $j > i$}
            \State $A_{\text{chordal}}[i,j] \gets 1$
            \State $A_{\text{chordal}}[j,i] \gets 1$
        \EndFor
    \EndFor
\EndFor

\State \textcolor{gray}{// Remove self-loops}
\For{$v = 1$ to $n$}
    \State $A_{\text{chordal}}[v,v] \gets 0$
\EndFor

\State \textcolor{gray}{// Reverse permutation to restore original ordering}
\State $\pi^{-1} \gets \text{inverse permutation of } \pi$
\State $A_{\text{chordal}} \gets A_{\text{chordal}}[\pi^{-1}, \pi^{-1}]$

\State \Return $A_{\text{chordal}}$
\end{algorithmic}
\end{algorithm}

\subsection{PSD Completion Algorithm}
\label{app:psd_completion}

The PSD completion algorithm reconstructs the full voltage product matrix $\mathbf{W}_t$ from the sparse clique-indexed variables $\{\mathbf{W}_{C_i, t}\}_{i=1,\dots,k}$ obtained from solving the chordal SDP relaxation. This procedure exploits the chordal structure of the sparsity graph to efficiently fill in the missing entries while preserving the positive semidefinite property.

\subsubsection*{Algorithm Overview}

Given a chordal graph $G = (\mathcal{V}, \mathcal{E})$ with maximal cliques $\mathcal{C} = \{C_1, \dots, C_k\}$ and a perfect elimination ordering (PEO) $\texttt{peo} = [v_1, \dots, v_n]$, the completion proceeds by processing vertices in reverse PEO order. At each step, we leverage the clique separator property: for any vertex $v$ being processed, its previously processed neighbors $U$ form a clique that \emph{separates} $v$ from all other previously processed vertices $T$, meaning that every path from $v$ to any vertex in $T$ passes through $U$ \cite{lauritzen1996, golumbic2004}. (Equivalently, removing $U$ from the graph disconnects $v$ from $T$.) This structure allows us to compute the missing entries between $v$ and $T$ using the Schur complement relationship \cite{grone1984, andersen2015}.

The key fact is that in a positive definite matrix consistent with a chordal sparsity pattern, if a clique $U$ separates vertex sets $S$ and $T$, then the matrix entries between $S$ and $T$ are completely determined by their interactions with $U$ \cite{grone1984, fukuda2001}. (In the language of Gaussian graphical models, this corresponds to the fact that variables in $S$ and $T$ are conditionally independent given $U$, so their covariance is determined entirely by the blocks involving $U$ \cite{lauritzen1996, WainwrightJordan2008}.) Formally, if we partition the matrix $\mathbf{W}_t$ as
\begin{align}
\mathbf{W}_t = 
\begin{pmatrix}
\mathbf{W}_{SS} & \mathbf{W}_{SU} & \mathbf{W}_{ST} \\
\mathbf{W}_{US} & \mathbf{W}_{UU} & \mathbf{W}_{UT} \\
\mathbf{W}_{TS} & \mathbf{W}_{TU} & \mathbf{W}_{TT}
\end{pmatrix},
\end{align}
where $S = \{v\}$, $U$ are the neighbors of $v$ among the processed vertices, and $T$ are the remaining processed vertices, then the missing blocks satisfy
\begin{align}
\mathbf{W}_{ST} &= \mathbf{W}_{SU} \mathbf{W}_{UU}^{-1} \mathbf{W}_{UT}, \label{completion_ST}\\
\mathbf{W}_{TS} &= \mathbf{W}_{TU} \mathbf{W}_{UU}^{-1} \mathbf{W}_{US}. \label{completion_TS}
\end{align}
This follows directly from the Schur complement identity and the clique separator property of chordal graphs \cite{andersen2015, Zhang2006}. Under these conditions, the completion is unique and preserves positive definiteness \cite{grone1984}.

\begin{algorithm}[H]
\caption{PSD Completion for Chordal Graphs}
\label{alg:psd_completion}
\begin{algorithmic}[1]
\Require Chordal graph $G = (\mathcal{V}, \mathcal{E})$, PEO $\texttt{peo} = [v_1, \dots, v_n]$, clique matrices $\{\mathbf{W}_{C_i, t}\}_{i=1,\dots,k}$
\Ensure Completed global matrix $\mathbf{W}_t \in \mathbb{C}^{n \times n}$

\State Initialize $\mathbf{W}_t \gets \mathbf{0}_{n \times n}$
\State \textcolor{gray}{// First pass: populate entries from clique constraints}
\For{each maximal clique $C_i \in \mathcal{C}$}
    \For{all $(v,w) \in C_i \times C_i$}
        \State $\mathbf{W}_t[v,w] \gets \mathbf{W}_{C_i,t}[v,w]$ \textcolor{gray}{// Copy clique entries to global matrix}
    \EndFor
\EndFor

\State Initialize set of processed vertices $U \gets \emptyset$

\State \textcolor{gray}{// Second pass: complete missing entries in reverse PEO order}
\For{$i = n-1$ down to $0$}
    \State $v \gets \texttt{peo}[i]$ \textcolor{gray}{// Current vertex being processed}
    \State $S \gets \{v\}$
    
    \State $U_v \gets \{u \in U \mid (u,v) \in \mathcal{E}\}$ \textcolor{gray}{// Neighbors of $v$ in $U$}
    \State $T_v \gets U \setminus U_v$ \textcolor{gray}{// Non-neighbors of $v$ in $U$}
    
    \If{$|T_v| > 0$ and $|U_v| > 0$} \textcolor{gray}{// Only fill if there are missing entries}
        \State \textcolor{gray}{// Extract relevant submatrices}
        \State $\mathbf{W}_{SU_v} \gets \mathbf{W}_t[S, U_v]$ \textcolor{gray}{// $1 \times |U_v|$ matrix}
        \State $\mathbf{W}_{U_vU_v} \gets \mathbf{W}_t[U_v, U_v]$ \textcolor{gray}{// $|U_v| \times |U_v|$ matrix}
        \State $\mathbf{W}_{T_vU_v} \gets \mathbf{W}_t[T_v, U_v]$ \textcolor{gray}{// $|T_v| \times |U_v|$ matrix}

        \State $\mathbf{W}_{U_vU_v}^\dagger \gets \text{pinv}(\mathbf{W}_{U_vU_v})$
        
        \State \textcolor{gray}{// Fill missing entries using Schur complement}
        \State $\mathbf{W}_t[S, T_v] \gets \mathbf{W}_{SU_v} \cdot \mathbf{W}_{U_vU_v}^\dagger \cdot \mathbf{W}_{T_vU_v}^\top$
        \State $\mathbf{W}_t[T_v, S] \gets \mathbf{W}_{T_vU_v} \cdot \mathbf{W}_{U_vU_v}^\dagger \cdot \mathbf{W}_{SU_v}^\top$
    \EndIf
    
    \State $U \gets U \cup \{v\}$ \textcolor{gray}{// Mark vertex as processed}
\EndFor

\State \Return $\mathbf{W}_t$
\end{algorithmic}
\end{algorithm}

\subsubsection*{Correctness and Complexity}

The correctness of this algorithm follows from the fundamental property of chordal graphs: in a perfect elimination ordering, each vertex and its unprocessed neighbors form a clique \citep{grone1984}. This ensures that the separating set $U_v$ in each iteration is indeed a clique, which guarantees that the Schur complement formulas \eqref{completion_ST}--\eqref{completion_TS} correctly reconstruct the PSD completion.

Let $\omega$ denote the treewidth of the chordal graph (i.e., the size of the largest clique minus one). 
If, when processing a vertex, the separating clique $U$ has size $k = |U|$, then the computation of the corresponding matrix blocks involves matrix operations on $O(k) \times O(k)$ dense submatrices, each requiring $\Theta(k^3)$ arithmetic operations. 
Hence, the total arithmetic cost of the completion procedure is
\[
\Theta\!\Big(\sum_{v \in \mathcal{V}} k_v^3\Big),
\]
where $k_v$ is the size of the separating clique associated with vertex $v$. 
Using the uniform bound $k_v \le \omega$ gives the worst-case complexity
\[
O(n\,\omega^3).
\]
For dense graphs, $\omega = \Theta(n)$, leading to $O(n^4)$ complexity, whereas for sparse graphs with small treewidth, the procedure is highly efficient. 
In power networks, the treewidth is typically small (often $\omega \le 10$), making the algorithm computationally tractable even for large-scale systems.

\end{document}